\documentclass[12pt]{article}   

\usepackage{graphicx}  
\usepackage{amssymb}
\usepackage{amsmath}
\usepackage{subcaption}
\usepackage{undertilde}
\usepackage{tabularx, makecell, multirow} 

\usepackage{rotating}
\usepackage{tikz}
\usetikzlibrary{arrows,
                positioning,
                shapes.geometric}

\usepackage{adjustbox}
                
\usepackage{hyperref}
\hypersetup{
    colorlinks=true,
    linkcolor=blue,
    citecolor=blue,
    filecolor= magenta,      
    urlcolor= magenta 
}

\newtheorem{theorem}{Theorem}[section]

\newcommand{\rr}{{\mathbf r}}
\newcommand{\s}{{\mathbf s}}

\newcommand{\nn}{{\mathbf n}}
\newcommand{\rp}{{\mathbf r}^{\prime}}

\newcommand{\Phit}{\tilde{\Phi}}

\newcommand{\R}{\mathbb{R}}

\newcommand{\ei}{\epsilon_\infty}
\newcommand{\es}{\epsilon_s}
\newcommand{\ep}{\epsilon_p}
\newcommand{\ez}{\epsilon_0}

\newcommand{\Pphi}[1]{\partial \Phi#1}

\title{A Nonlocal size modified Poisson-Boltzmann Model and Its Finite Element Solver for \\
Protein in Multi-Species Ionic Solution}
\author{Dexuan Xie\thanks{Corresponding author: Dexuan Xie (dxie@uwm.edu),
Department of Mathematical Sciences, University of Wisconsin-Milwaukee, Milwaukee, WI, 53201-0413, USA}, Liam Jemison\thanks{Department of Mathematical Sciences, University of Wisconsin-Milwaukee, Milwaukee, WI, 53201-0413, USA}, and Yi Jiang\thanks{Shenzhen Institute of Advanced Technology, Chinese Academy of Sciences, Shenzhen, China.}}

\normalsize

\date{}

\begin{document}						

\maketitle

\begin{abstract}
The Poisson-Boltzmann (PB) model is a widely used implicit solvent model in protein simulations. Although variants, such as the size modified PB and nonlocal modified PB models, have been developed to account for ionic size effects and nonlocal dielectric correlations, no existing PB variants simultaneously incorporate both, due to significant modeling and computational challenges. To address this gap, in this paper, a nonlocal size modified Poisson-Boltzmann (NSMPB) model is introduced and solved using a finite element method for a protein with a three-dimensional molecular structure and an ionic solution containing multiple ion species. In particular, a novel solution decomposition is proposed to overcome the difficulties caused by the increased nonlinearity, nonlocality, and solution singularities of the model. It is then applied to the development of the NSMPB finite element solver, which includes an efficient modified Newton iterative method, an effective damping parameter selection strategy, and good selections of initial iterations. Moreover, the construction of the modified Newton iterative method is mathematically justified. Furthermore, an NSMPB finite element package is developed by integrating a mesh generation tool, a protein data bank file retrieval program, and the PDB2PQR package to simplify and accelerate its usage and application. Finally, numerical experiments are conducted on an ionic solution with four species, proteins with up to 11,439 atoms, and irregular interface-fitted tetrahedral box meshes with up to 1,188,840 vertices. The numerical results confirm the fast convergence and strong robustness of the modified Newton iterative method, demonstrate the high performance of the package, and highlight the crucial roles played by the damping parameter and initial iteration selections in enhancing the method’s convergence. The package will be a valuable tool in protein simulations.
\end{abstract}



\section{Introduction}
One fundamental challenge in protein simulation is the development of ionic solvent models, as proteins exist naturally in ionic solutions. The explicit solvent approach is widely used in molecular dynamics (MD) simulations, facilitated by popular software packages such as CHARMM \cite{charmm}, AMBER \cite{cornell95}, GROMACS \cite{gromacs2008} and NAMD \cite{namd96}. However, this approach models water molecules as explicit variables in the MD system, leading to high computational costs. To address this, the implicit solvent approach, also known as dielectric continuum modeling, was introduced \cite{honig95,roux99}. By treating water as a dielectric medium, this approach significantly reduces the complexity of protein simulations. Among the various implicit solvent models, the Poisson-Boltzmann (PB) model is the most commonly used. It treats the water solvent and protein regions as two distinct continuum dielectrics with different permittivity constants \cite{holst2001,xiePBE2013}. The PB model has been extensively applied to computing electrostatic solvation and binding free energies, studying biomolecular electrostatics, and advancing bioengineering applications such as protein docking and rational drug design \cite{fogolari2002poisson,McCammonPBEreivew2008,Vizcarra2005}.

To enhance the accuracy of PB models in computing electrostatic solvation free energies, size modified PB (SMPB) models have been developed to account for ionic size effects in solutions. Along these lines, we have contributed to the development of SMPB models and their efficient finite element solvers \cite{xiePBE2013, SDPBS2015, Yi_Xie2014,Li-Xie2014, XieLi2014, Ying-Xie2014,nuSMPBE2017}. Models of ionic solutions that treat ions as spheres exhibit entropy and energy characteristics significantly different from those that consider ions as point charges. For example, the ability of nerve cells to distinguish between sodium ions (Na$^+$) and potassium ions (K$^+$) -- despite their identical charges --- depends critically on their distinct sizes. This underscores the importance of size differentiation in dielectric continuum models.  

Moreover, the dielectric properties of water solvents are inherently nonlocal, arising from hydrogen bonding among water molecules, which induces strong polarization correlations, particularly near protein surfaces where ion concentrations are highest. This phenomenon has driven the development of nonlocal dielectric models \cite{PhysRevLett.79.3435, kornyshev1978model, vorotyntsev1978model,Basilevsky2008}. However, nonlocal dielectric continuum models pose significant computational challenges due to their governing integro-differential equations, which involve convolutions of the gradient vector of the electrostatic potential over three-dimensional space. Early studies mainly addressed simplified cases, such as single ions immersed in water \cite{PhysRevLett.79.3435,Basilevsky2008,ref:nonlocelectrostatcavity,ref:nonlocontindielsuscept,ref:nonlocontmodnonrigicav,NonlocalTheory1,dielectric1985,NonlocalTheory6,rubinstein2004influence,lrsBIBfi}. Recent advances have substantially alleviated computational burden, extending their application from water solvents to ionic solutions, and from individual ions to proteins immersed in ionic solutions \cite{PhysRevLett.93.108104,Weggler20104059,xie2011nonlocal,xie_volkmer2011,xie-nonlocal-solver2012,HansXie2015b, ScottXie2015}. These developments have led to a nonlocal modified Poisson-Boltzmann (NMPB) model and its finite element solvers for a protein in a symmetric 1:1 ionic solution \cite{ying2018accelerated,xie-nonlocal2014}. Comparisons with local models, using experimental data and analytical solutions, have confirmed the superior predictive power of nonlocal dielectric models \cite{HansXie2015b,ying2018accelerated}.  

Despite these advancements, no existing Poisson-Boltzmann (PB) variant has incorporated both ionic-size effects and nonlocal dielectric properties due to inherent modeling and computational challenges. To bridge this gap, we propose a nonlocal size modified Poisson-Boltzmann (NSMPB) model in this work and develop a finite element iterative method to solve it. This method is applicable to a protein with a three-dimensional molecular structure immersed in an ionic solution containing multiple ion species within a box domain. As a special case, setting all ionic sizes to zero reduces the NSMPB model to the NMPB model reported in \cite{xie-nonlocal2014}, producing an improved NMPB finite element solver that works for both symmetric 1:1 ionic solutions and mixture solutions with multiple ionic species.

The NSMPB model turns out to be significantly more challenging to solve numerically than either the SMPB model or the NMPB model due to its stronger nonlinearity and more complex nonlocal terms, in addition to the solution singularity induced by the atomic charges of the protein and discontinuous interface conditions. Thus, new mathematical and numerical techniques are required to develop an effective NSMPB finite element solver. 

To do so, we begin by decomposing the NSMPB solution $u$ into the sum of three component functions, denoted by $G$, $\Psi$, and $\Phit$ as done in our previous work. Here, $G$ is a known function that collects all singular points of $u$. This decomposition allows us to construct a linear nonlocal interface boundary value problem for computing $\Psi$ and a nonlinear nonlocal interface boundary value problem for computing $\Phit$, neither of which involves singularity problems. To this end, we can entirely circumvent the singularity problems induced by atomic charges.

We next focus on the development of a modified Newton iterative method to solve the nonlinear problem for $\Phit$, since we can adopt the finite element method reported in \cite[Eq. 33]{xie-nonlocal2014} to numerically solve the linear problem for $\Psi$. By treating the convolution of $\Phit$ as an unknown function, we novelly formulate the nonlinear problem into a nonlinear finite element variational system in terms of $\Phit$ and its convolution, avoiding the numerical difficulty of computing the convolution of $\Phit$. Unlike the NMPB case \cite[Eq. 40]{xie-nonlocal2014}, this nonlinear variational system presents significantly greater challenges in constructing a modified Newton iterative method. We overcome these difficulties to derive an effective and efficient modified Newton iterative method, complemented by a simple yet effective damping parameter selection scheme, outlined in Algorithm 1, and two well-chosen initial iterations, which can significantly enhance the convergence of the method. We provide a detailed construction of the method and rigorously justify it in Theorems 3.2 to 3.4. In addition, we consider two more initial iteration selections to demonstrate the robustness of our modified Newton iterative method. Furthermore, we linearize the nonlinear problem to derive a linear NSMPB model, which is valuable in its own right due to its computational simplicity.

We summarize our NSMPB finite element solver in Algorithm~2 and implement it as a program package based on the finite element library of the FEniCS project (Version 2019.1.0) \cite{dolfin2012}, our SMPB program package \cite{nuSMPBE2017}, and our NMPB program package \cite{xie-nonlocal-solver2012}. To ensure modularity and reusability, we design the NSMPB program package using object-oriented programming techniques. This allows us to seamlessly integrate our NSMPB finite element program with a mesh generation package, a Python program for downloading PDB files from the Protein Data Bank (PDB, {\em http://www.rcsb.org/}) using a four-character PDB identification code, and the PDB2PQR package ({\em http://www.poissonboltzmann.org/pdb2pqr/}) for converting PDB files to PQR files \cite{dolinsky2004pdb2pqr}. A PQR file supplements missing data from the original PDB file, such as hydrogen atoms, atomic charge numbers, and atomic radii, which are required by the mesh generation package. To further improve computational efficiency, we wrote Fortran subroutines to calculate computationally intensive terms and functions (such as convolutions and nonlinear system residual vectors) and converted them into Python modules using the Fortran-to-Python interface generator \texttt{f2py} ({\em https://numpy.org/doc/stable/f2py/index.html}) to call them directly by our NSMPB package. Through these efforts, we create a user-friendly NSMPB finite element package, where users only need to input a protein PDB ID to generate a finite element solution of the NSMPB model (an electrostatic potential function) and its convolution and ionic concentration functions. By automating these processes, we significantly reduce the technical barriers that we have during the development of the package. Consequently, the package will provide researchers with a powerful tool for computing electrostatic solvation free energies, especially in the cases where the standard mean-field approach falls short because of its neglect of ionic size effects and nonlocal dielectric properties.

Finally, we conducted numerical experiments using the NSMPB finite element package for three proteins with up to 11,439 atoms in a mixture solution containing four ionic species. We used six irregular interface-fitted box domain meshes with up to 1,188,840 vertices and 7,278,073 tetrahedra. The numerical results highlight the high performance of the NSMPB finite element package, the fast convergence rate of our modified Newton iterative method, and the critical role of our damping parameter and initial iteration selections in enhancing the method's convergence. For example, our modified Newton iterative method required only 9 iterations to reduce the absolute residual error of a nonlinear finite element system from approximately $10^2$ to $10^{-7}$, as detailed in Figure~\ref{convergence4mesh}, in just 64 seconds on an M4 chip of Mac mini, as given in Table~2, for a protein (PDB ID: 1CID) on a box mesh with 193,592 vertices (case of Mesh~5). In this test, the nonlinear finite element system has about 387,184 unknowns since it has two unknown functions, each having 193,592 unknowns. These numerical results confirm the fast convergence rate of our modified Newton iterative method and demonstrate the high performance of our NSMPB finite element program package.

The remaining sections of the paper are outlined as follows.  In Section 2, we present the NSMPB model. In Section 3, we present the NSMPB finite element iterative method. In Section 4, we report the NSMPB program package and numerical results. Conclusions are made in Section 5.

\section{A nonlocal size modified Poisson-Boltzmann model}
\begin{figure}[t]
\centering
 \includegraphics[width=0.3\textwidth]{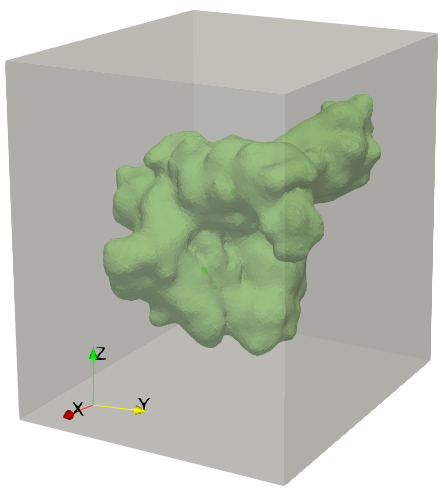}
\caption{An illustration of the box domain $\Omega$ partitioned into the protein region $D_p$ (green), the solvent region $D_s$ (gray), and the interface $\Gamma$ between $D_p$ and $D_s$. In this figure, $\Gamma$ is set as a molecular surface of the protein with the protein data bank identifier (PDB ID) 1C4K.      }       
\label{fig-box_partition}          
\end{figure}

We select a sufficiently large box domain, $\Omega$, satisfying the domain partition:
\[ \Omega = D_{p} \cup D_{s} \cup \Gamma, \]
where $D_{p}$, $D_{s}$, and $\Gamma$ denote a protein region, a solvent region, and an interface between $D_{p}$ and $D_{s}$, respectively. In particular, $D_p$ is surrounded by $D_s$, $D_s$ contains an ionic solution with $n$ different ionic species, and $\Gamma$ is a molecular surface of the protein, which wraps a three-dimensional protein structure with $n_p$ atoms. The interface $\Gamma$ can also be treated as the boundary of $D_p$.  An illustration of the domain partition is given in Figure~\ref{fig-box_partition}, showing that $D_p$, $D_s$, and $\Gamma$ can have complicated geometric shapes.

Let $\Phi$ denote the electrostatic potential function of the electric field induced by the atomic charges in $D_p$ and ionic charges in $D_s$, and let $c_i$ denote the ionic concentration function of the $i$th species for $i=1,2,\ldots,n$. Note that $\Phi$ and $c_i$ are defined in two different domains, $\Omega$ and $D_s$, respectively. In the international system of units, $\Phi$ has units of volts (V) and $c_i$ is in the number of ions per meter cubed (1/m$^3$) since length is measured in meters (m). When the atomic charge number $z_j$ and atomic position $\rr_j$ of atom $j$ and the charge number $Z_i$ of ionic species $i$ are given, from the nonlocal continuum implicit solvent theory \cite{xie-nonlocal-solver2012}, we can obtain a nonlocal Poisson dielectric model as follows:
\begin{equation}
\label{phi-u1}
  \left\{\begin{array}{cl}
    -\epsilon_p\Delta\Phi(\rr) = \frac{e_{c}}{\epsilon_0}\sum\limits_{j=1}^{n_{p}}z_{j} \delta_{\rr_{j}}, & \rr \in D_p,\\
-\epsilon_{\infty}\Delta\Phi(\rr) + \frac{\epsilon_s-\epsilon_{\infty}}{\lambda^2}[\Phi(\rr)-(\Phi\ast Q_{\lambda})(\rr)] 
= \frac{e_{c}}{\epsilon_0}\sum\limits_{i=1}^{n}Z_{i}c_{i}(\rr), & \rr \in D_s,\\
\Phi(\s^{-}) = \Phi(\s^{+}), \quad  \epsilon_p\frac{\partial\Phi(\s^{-})}{\partial \nn(\s)} = \epsilon_\infty\frac{\partial\Phi(\s^{+})}{\partial \nn(\s)}+(\epsilon_s-\epsilon_\infty)\frac{\partial (\Phi\ast Q_{\lambda})(\s)}{\partial \nn(\s)}, & \s\in\Gamma,\\
   \Phi(\s) = \bar{g}(\s), & \s \in \partial \Omega,
  \end{array}
\right.
\end{equation}
where $D_{p}$ and $D_{s}$ have been treated as dielectric media with two different relative dielectric constants $\ep$ and $\es$, respectively; $\lambda$ is a parameter for characterizing the polarization correlations of water molecules or the spatial-frequency dependence of a dielectric medium in a more general sense; $\ei$ is a relative dielectric constant satisfying $\ei < \es$, which corresponds to the case $\lambda\to \infty$; $e_c$ is the elementary charge in Coulomb (C); $\ez$ is the permittivity of the vacuum in Farad/meter (F/m); $\bar{g}$ is a boundary value function; $\partial \Omega$ denotes the boundary of $\Omega$; $\delta_{\rr_{j}}$ is the Dirac delta distribution at  $\rr_{j}$; $\nn$ denotes the unit outward normal direction of $D_p$; $\frac{\Pphi(\s)}{\partial \nn(\s)}=\nabla \Phi(\s)\cdot \nn(\s)$; and  $\Phi\ast Q_{\lambda}$ denotes the convolution of $\Phi$ with the kernel function $Q_{\lambda}$, which is defined by
\begin{equation}\label{Q-def}
  \Phi\ast Q_{\lambda} =  \int_{\mathbb{R}^3}Q_{\lambda}(\rr-\rp) \Phi(\rp)\mathrm{d}\rp
\qquad \mbox{with } \quad Q_{\lambda}(\rr)=\frac{e^{-|\rr|/\lambda}}{4\pi \lambda^2 |\rr|}. 
 \end{equation}
Here, $\mathbb{R}^3$ denotes the three-dimensional Euclidean space and $\rr$ represents the position vector of a point in the space. 

Note that in molecular simulation, length data are often given in angstroms (\AA), and  $c_{i}$ in  moles per liter (mol/L). Hence, for the convenience of calculation, we change the length unit from meters to angstroms (\AA) and the concentration unit  from mol/L   into 1/\AA$^3$ by the transformations:  
\begin{equation*}
  1 m = 10^{10} \text{\AA}, \qquad 1 L = m^3/10^3,  \qquad  
    \mbox{ mol /L} =  10^{3} N_{A} /m^{3}= 10^{-27} N_{A} /\mbox{\r{A}}^{3}, 
\end{equation*}
where $N_A$ is the  Avogadro number for estimating  the number of ions per mole (i.e., 1 mole = 1 $N_A$).
We then rescale the potential $\Phi$ to a dimensionless potential function, $u$, by the formula
\begin{equation}
\label{unit-changes}
   u= \frac{e_{c}}{k_{B}T}\Phi,
\end{equation}
where $k_B$ is the Boltzmann constant in Joule/Kelvin (J/K), and $T$ is the absolute temperature in Kelvin (K). Using the above unit changes and \eqref{unit-changes}, we can get an interface boundary value problem for defining $u$ as follows:
\begin{equation}
\label{nonlocal-Poisson-u}
 \left\{\begin{array}{cl}
     -\epsilon_p\Delta u(\rr) =\alpha \sum\limits_{j=1}^{n_{p}}z_{j}  \delta_{\rr_{j}}, & \rr \in D_p,\\
-\epsilon_{\infty}\Delta u + \frac{\epsilon_s-\epsilon_\infty}{\lambda^2}[u-(u\ast Q_{\lambda})]
= \beta \sum\limits_{i=1}^{n}Z_{i}c_i(\rr), & \rr \in D_s,\\
u(\s^{-}) = u(\s^{+}), \quad
\epsilon_p\frac{\partial u(\s^{-})}{\partial \nn(\s)} = \epsilon_\infty\frac{\partial u(\s^{+})}{\partial \nn(\s)}+(\epsilon_s-\epsilon_\infty)\frac{\partial (u\ast Q_{\lambda})(\s)}{\partial \nn(\s)}, & \s\in\Gamma,\\
   u(\rr) = g(\s), & \s \in \partial \Omega,
\end{array}
\right.
\end{equation}  
where $\Delta$ denotes the Laplace operator, $\alpha$ and $\beta$ are the two constants defined by
 \begin{equation}
\label{alpha-beta}
\alpha = \frac{10^{10}e_{c}^{2}}{\ez k_{B}T}, \qquad \beta = \frac{N_A e_{c}^{2}}{10^{17}\ez k_{B}T},
 \end{equation}  
$\nn$ is the unit outward normal vector of $D_{p}$, $\frac{\partial u(\s)}{\partial \nn(\s)}=\nabla u(\s)\cdot \nn(\s)$ with $\nabla$ denoting the gradient operator, $g = \frac{e_{c}}{k_{B}T} \bar{g}$, $u(\s^{\pm}) = \lim_{t\rightarrow 0^+} u(\s\pm t\nn(\s))$, and
$\frac{\partial u(\s^{\pm})}{\partial \nn(\s)} = \lim_{t\rightarrow 0^+} \frac{\partial u(\s\pm t\nn(\s)))}{\partial \nn(\s)}$.

For $T=298.15$ K, $k_B = 1.380648813\times10^{-23}$ J/K, 
$N_A=6.02214129\times10^{23}$, $e_c=1.602176565\times10^{-19} \mbox{ C}$,  and  
$\ez = 8.854187817\times10^{-12} \mbox{ F/m}$,    the values of $\alpha$ and $\beta$ can be estimated as follows:
\[ \alpha\approx 7042.93990033, \qquad  \beta  \approx 4.24135792.
  \]

Different selections of ionic concentration functions may yield different electrostatic potential functions. To reflect ionic size effects, a selection of concentrations $c_i$ for $i=1,2,\ldots,n$ has been made in \cite[Eq. 27]{nuSMPBEic} as follows:
\begin{equation}
\label{concentration4nsmpb}
   c_{i}  = \frac{ c_{i}^b e^{-Z_{i} u(\rr)}}
         {1 + \gamma  \frac{\bar{v}^2}{v_0} \sum\limits_{j=1}^{n} 
         c_{j}^{b}e^{-Z_{j} u}},   \quad i=1,2,\ldots,n,
\end{equation}
where $\bar{v} = \frac{1}{n}\sum_{j=1}^n v_j$ with $v_i$ denoting the volume of an  ion of species $i$, $\gamma = 10^{-27} N_A$, which  is about $6.02214129\times10^{-4}$, $v_0$ is a size scaling parameter (e.g.,  $v_0 = \min_{1\leq i\leq n} v_i$ by default), and $c_i^b$ denotes a bulk concentration of ionic species $i$ in moles per liter.

Substituting \eqref{concentration4nsmpb} to the nonlocal Poisson dielectric model \eqref{nonlocal-Poisson-u}, we obtain a nonlocal size modified Poisson-Boltzmann (NSMPB) model as follows:
\begin{equation}
\label{NSMPB-def}
 \left\{\begin{array}{cl}
     -\epsilon_p\Delta u(\rr) =\alpha \sum\limits_{j=1}^{n_{p}}z_{j}  \delta_{\rr_{j}}, & \rr \in D_p,\\
\displaystyle -\epsilon_{\infty}\Delta u + \frac{\epsilon_s-\epsilon_\infty}{\lambda^2}[u-(u\ast Q_{\lambda})]
- \beta \frac{ \sum\limits_{i=1}^{n}Z_{i} c_{i}^{b} e^{-Z_{i} u(\rr)}}
{1+ \gamma  \frac{\bar{v}^2}{v_0}  \sum\limits_{i=1}^n c_{i}^{b}e^{-Z_{i} u(\rr)}} =0, & \rr \in D_s,\\
\displaystyle u(\s^{-}) = u(\s^{+}), \quad
\epsilon_p\frac{\partial u(\s^{-})}{\partial \nn(\s)} = \epsilon_\infty\frac{\partial u(\s^{+})}{\partial \nn(\s)}+(\epsilon_s-\epsilon_\infty)\frac{\partial (u\ast Q_{\lambda})(\s)}{\partial \nn(\s)}, & \s\in\Gamma,\\
   u (\s) =  g(\s), & \s \in \partial \Omega.
\end{array}
\right.
\end{equation}  

Setting all ion sizes $v_j=0$ for $j=1,2, \ldots, n$, we reduce the NSMPB model to the NMPB model, reported in \cite[Eq. (21)]{xie-nonlocal2014}, as follows:
\begin{equation}
\label{NMPB-def}
 \left\{\begin{array}{cl}
     -\epsilon_p\Delta u(\rr) =\alpha \sum\limits_{j=1}^{n_{p}}z_{j}  \delta_{\rr_{j}}, & \rr \in D_p,\\
\displaystyle -\epsilon_{\infty}\Delta u + \frac{\epsilon_s-\epsilon_\infty}{\lambda^2}[u-(u\ast Q_{\lambda})]
- \beta \sum\limits_{i=1}^{n}Z_{i} c_{i}^{b} e^{-Z_{i} u(\rr)} =0, & \rr \in D_s,\\
\displaystyle u(\s^{-}) = u(\s^{+}), \quad
\epsilon_p\frac{\partial u(\s^{-})}{\partial \nn(\s)} = \epsilon_\infty\frac{\partial u(\s^{+})}{\partial \nn(\s)}+(\epsilon_s-\epsilon_\infty)\frac{\partial (u\ast Q_{\lambda})(\s)}{\partial \nn(\s)}, & \s\in\Gamma,\\
   u (\s) =  g(\s), & \s \in \partial \Omega.
\end{array}
\right.
\end{equation}  
In this sense, the NSMPB model can be regarded as an extension of the NMPB model to reflect the effects of ionic size in the calculation of ionic concentrations. 

The problem now becomes to solve the interface boundary value problem for $u$ since we can obtain $n$ ionic concentration functions, $c_i$ for $i=1,2,\ldots,n$, using the formula \eqref{concentration4nsmpb} when $u$ is given.

{\bf Remark:}
Setting $\epsilon_{\infty}=\es$ (i.e. without considering any nonlocal effects), we can reduce the NSMPB model \eqref{NSMPB-def} to a size modified PB (SMPB) model, reported in \cite[Eq. (36)]{nuSMPBE2017}, as follows:
\begin{equation}
\label{SMPB-def}
 \left\{\begin{array}{cl}
     -\epsilon_p\Delta u(\rr) =\alpha \sum\limits_{j=1}^{n_{p}}z_{j}  \delta_{\rr_{j}}, & \rr \in D_p,\\
\displaystyle \epsilon_{s}\Delta u +
\beta \frac{ \sum\limits_{i=1}^{n}Z_{i} c_{i}^{b} e^{-Z_{i} u(\rr)}}
{1+ \gamma  \frac{\bar{v}^2}{v_0}  \sum\limits_{i=1}^n c_{i}^{b}e^{-Z_{i} u(\rr)}} =0, & \rr \in D_s,\\
\displaystyle u(\s^{-}) = u(\s^{+}), \quad
\epsilon_p\frac{\partial u(\s^{-})}{\partial \nn(\s)} = \epsilon_\infty\frac{\partial u(\s^{+})}{\partial \nn(\s)}, & \s\in\Gamma,\\
   u (\s) =  g(\s), & \s \in \partial \Omega.
\end{array}
\right.
\end{equation}  
Thus, the SMPB model can be regarded as a special case of the NSMPB model. 

\section{An NSMPB finite element iterative method} 

In this section, we present a finite element iterative method for solving the NSMPB model \eqref{NSMPB-def}. To overcome the singularity difficulty caused by atomic charges, we start with a decomposition of $u$ by 
\begin{equation}\label{SplitForm-u}
u(\rr) = G(\rr) +  \Psi(\rr) + \Phit(\rr), \quad \rr \in \Omega,
\end{equation}
where  $G$ is given by the algebraic expression
\begin{equation}
\label{Def_G}
   G(\rr)=\frac{\alpha}{4\pi\ep}  \sum\limits_{j=1}^{n_p}\frac{z_{j}}{| \rr-\rr_j |},
\end{equation}
 $ \Psi$ is a solution of the linear nonlocal interface boundary value problem
 \begin{equation}\label{Psi_b_nonlocal}
\left\{
\begin{array}{cl}
\Delta\Psi(\rr) = 0, & \rr \in D_p,\\
-\epsilon_\infty\Delta\Psi(\rr) +  \frac{\epsilon_s-\epsilon_\infty}{\lambda^2}[\Psi(\rr)-(\Psi \ast Q_{\lambda})(\rr)]= - \frac{\epsilon_s-\epsilon_\infty}{\lambda^2}[G(\rr)-(G\ast Q_{\lambda})(\rr)], & \rr \in D_s,\\
\Psi(\s^{-}) =\Psi(\s^{+}),  \quad 
  \epsilon_p\frac{\partial\Psi(\s^{-})}{\partial \nn(\s)} -  \epsilon_{\infty}\frac{\partial\Psi(\s^{+})}{\partial \nn(\s)} 
 =  (\epsilon_s-\epsilon_\infty)\frac{\partial (\Psi \ast Q_{\lambda})(\s)}{\partial \nn(\s)} +g_{\Gamma}(\s), &  \s\in\Gamma,\\
\Psi(\s) = g(\s) - G(\s), & \s \in \partial \Omega,
\end{array}
\right.
\end{equation}
and $\Phit$ is a solution of the nonlinear nonlocal interface boundary value problem:
\begin{equation}
\label{Phit-NSMPB}
\left\{\begin{array}{cl}
\Delta \Phit (\rr) = 0, & \rr\in D_p,\\
-\epsilon_\infty\Delta \Phit(\rr) + \frac{\epsilon_s-\epsilon_\infty}{\lambda^2}[\Phit(\rr) - (\Phit  \ast Q_{\lambda})(\rr)] -
\beta \frac{ \sum\limits_{i=1}^{n}Z_{i} c_{i}^{b}w_i(\rr) e^{-Z_{i}\Phit(\rr)}}
{1+ \gamma  \frac{\bar{v}^2}{v_0}  \sum\limits_{i=1}^n c_{i}^{b}w_i(\rr)e^{-Z_{i}\Phit(\rr)}} =0, & \rr\in D_s,\\
 \Phit (\s^{-}) = \Phit (\s^{+}), \quad
\epsilon_p\frac{\partial \Phit (\s^{-})}{\partial \nn(\s)} - \epsilon_{\infty}\frac{\partial \Phit (\s^{+})}{\partial \nn(\s)} = (\epsilon_s-\epsilon_\infty)\frac{\partial ( \Phit  \ast Q_{\lambda})(\s)}{\partial \nn(\s)}, & \s\in\Gamma,\\
\Phit (\s) = 0, & \s \in \partial \Omega.
\end{array}\right.
\end{equation}
Here, $w_i(\rr)$ and $g_{\Gamma}$ are defined by
\begin{equation}
\label{g-def}
w_i(\rr) = e^{-Z_{i}[G(\rr) + \Psi(\rr)]}, 
\quad
g_{\Gamma}(\s) =  (\epsilon_s-\epsilon_\infty)\frac{\partial (G\ast Q_{\lambda})(\s)}{\partial \nn(\s)} + (\epsilon_\infty-\epsilon_p)\frac{\partial G(\s)}{\partial \nn(\s)},
\end{equation}
and $\frac{\partial G(\s)}{\partial \nn(\s)} = \nabla G(\s) \cdot \nn(\s)$ with $\nabla G(\s)$ being given by
\begin{equation}
\label{nable_G-def}
   \nabla G(\s) = -\frac{\alpha}{4\pi\ep} \sum_{j=1}^{n_p}z_{j}\frac{(\s-\rr_j)}{|\s-\rr_j|^3}.
\end{equation}
For brevity, the convolution $G\ast Q_{\lambda}$ is denoted by $\hat{G}$. The algebraic expressions of $\hat{G}$ and $\nabla \hat{G}$ can be found as follows:
\begin{equation}\label{Def-hatG}
   \hat{G}(\rr) = \frac{\alpha}{4\pi \ep}\sum^{n_p}_{j=1}z_j\frac{1-e^{-\frac{|\rr-\rr_j|}{\lambda}}}{|\rr-\rr_j|}, \quad 
    \nabla \hat{G}(\rr) = \frac{\alpha}{4\pi \ep}\sum^{n_p}_{j=1}z_j\frac{\left(1 + \frac{| \rr - \rr_j|}{\lambda} \right)e^{-\frac{|\rr-\rr_j|}{\lambda}} - 1}{|\rr-\rr_j|^3}(\rr - \rr_j).
\end{equation}
Note that $G$ has collected all the singular points of $u$. Hence, the interface boundary value problems \eqref{Psi_b_nonlocal} and \eqref{Phit-NSMPB} do not involve any singularity points of $u$. As a result, they can be solved numerically much more easily than the original NSMPB model \eqref{NSMPB-def}. Consequently, the complexity of solving the original NSMPB model \eqref{NSMPB-def} has been markedly reduced due to the solution decomposition formula \eqref{SplitForm-u}.

We also note that $\Psi$ is independent of $\Phit$. Thus, we can calculate it before searching for $\Phit$. In this way, we can treat $\Psi$ as a known function when we solve problem \eqref{Phit-NSMPB} for $\Phit$. The problem \eqref{Psi_b_nonlocal} has been solved by adopting the linear finite element method reported in \cite[Eq. 33]{xie-nonlocal2014}. Hence, in this work, we focus on the construction of a finite element iterative method for solving \eqref{Phit-NSMPB}. After finding a finite element approximation of $\Phit$, we construct a finite element solution of the NSMPB model \eqref{NMPB-def} by \eqref{SplitForm-u}, and calculate each ionic concentration function using the expression \eqref{concentration4nsmpb}. 

For clarity, we present our finite element iterative method for solving \eqref{Phit-NSMPB} in the following five subsections. 

\subsection{Formulation of a nonlinear finite element variational problem}
To reformulate the problem \eqref{Phit-NSMPB} as a nonlinear finite element variational problem for computing $\Phit$ approximately, we first generate an interface-fitted tetrahedral mesh, $\Omega_h$, of $\Omega$. With the mesh, we construct a linear Lagrange finite element function space,  ${\cal M}$, such that ${\cal M} \subset H^{1}(\Omega)$. Here, $H^{1}(\Omega)$ denotes the usual Sobolev function space \cite{brennerScott}, and a function of ${\cal M}$ is linear in each tetrahedron of the mesh $\Omega_h$. We then define a subspace, ${\cal M}_{0}$, of  ${\cal M}$ by
\begin{equation}
\label{M0-def}
 {\cal M}_{0} = \{ v\in {\cal M} \; | \; v=0 \mbox{ on } \partial \Omega_h \},
\end{equation}
where $\partial \Omega_h$ is a triangular surface mesh of the boundary $\partial \Omega$ of box domain $\Omega$. By the reformulation techniques used in the construction of \cite[Eq. (5.4)]{xie-nonlocal-solver2012}, we can derive a finite element variational problem of \eqref{Phit-NSMPB}  as follows: 
\begin{equation}
\label{Phit-NSMPB-weak}
 \mbox{Find $\Phit  \in  {\cal M}_{0}$ such that }  \qquad  b(\Phit, v) =0 \qquad\forall v \in  {\cal M}_{0},
\end{equation}
where $b(\Phit,v)$ is a nonlinear functional of $\Phit$ given in the expression
\begin{equation}
\label{b-def}
\begin{aligned}
b(\Phit, v)  = &\ep \int_{D_{p}}\nabla \Phit (\rr)\cdot\nabla v(\rr)\mathrm{d}\rr +\ei \int_{D_{s}}\nabla \Phit (\rr)\cdot\nabla v(\rr)\mathrm{d}\rr   \\
+ & (\epsilon_s-\epsilon_\infty)\int_{D_s}\nabla (\Phit*Q_{\lambda})\cdot\nabla v \mathrm{d}\rr				    
    - \beta \int_{D_{s}}  \frac{ \sum\limits_{i=1}^{n}Z_{i} c_{i}^{b}w_i(\rr) e^{-Z_{i}\Phit(\rr)}}
{1+ \gamma  \frac{\bar{v}^2}{v_0}  \sum\limits_{i=1}^n c_{i}^{b}w_i(\rr)e^{-Z_{i}\Phit(\rr)}} v \mathrm{d}\rr.
\end{aligned}
\end{equation}

\subsection{A modified Newton iterative method}
One key step in developing a Newton iterative method to solve the nonlinear problem \eqref{Phit-NSMPB-weak} is to derive a linear expansion of the nonlinear functional $b(\Phit+p,v)$ in terms of $p$. Following the proof of Theorem~4.1 in \cite{xie-nonlocal2014}, we can get the linear expansion as done in the following theorem.

\begin{theorem}
\label{b_expansion}
Let $b(\Phit,v)$ be defined in \eqref{b-def} and $\| \cdot \|_{H^{1}(\Omega)}$ denote the norm of $H^{1}(\Omega)$. Then a linear expansion of $b(\Phit+p,v)$ is given as follows:
\begin{equation}
\label{b-expansion}
  b(\Phit+p,v) = b(\Phit,v) + b^{\prime}(\Phit,v;p) + O(\|p\|^{2}_{H^{1}(\Omega)}), \quad p \in {\cal M}_{0},
\end{equation}
where $b^{\prime}(\Phit,v;p)$ is the Fr\'echet derivative of the nonlinear functional $b(\Phit,v)$ at $p$, which is given in the expression 
\begin{equation}
\label{b-derivative}
\begin{aligned}
b^{\prime}(\Phit,v;p) & = \ep \int_{D_{p}}\nabla p(\rr)\cdot\nabla v(\rr)\mathrm{d}\rr +\ei \int_{D_{s}}\nabla p(\rr)\cdot\nabla v(\rr)\mathrm{d}\rr   \\
+ & (\epsilon_s-\epsilon_\infty)\int_{D_s}\nabla (p*Q_{\lambda})(\rr)\cdot\nabla v(\rr) \mathrm{d}\rr	\\	
+ & \beta \int_{D_s}  \frac{ A_1(\Phit) A_3(\Phit) -   \gamma  \frac{\bar{v}^2}{v_0}  \big[A_2(\Phit)\big]^2}{ \big[A_1(\Phit)\big]^2} p(\rr)v(\rr) d\rr.
\end{aligned}
\end{equation}
Here $A_1$, $A_2$, and $A_3$ are defined by
\[  A_1(\Phit) = 1+  \gamma  \frac{\bar{v}^2}{v_0} \sum\limits_{j=1}^n c_{j}^{b} w_i(\rr) e^{-Z_{j}\Phit(\rr)},\quad 
   A_2(\Phit)=\sum\limits_{i=1}^{n}Z_{i}  c_{i}^{b} w_i(\rr)e^{-Z_{i}\Phit(\rr)}, \]
and 
\[  A_3(\Phit) = \sum\limits_{i=1}^{n}Z_{i}^2  c_{i}^{b}w_i(\rr)e^{-Z_{i}\Phit(\rr)}.\]
\end{theorem}

We now construct the modified Newton iterative method using the linear expansion \eqref{b-expansion}.

Let $\Phit^{(0)}$ denote an initial guess to a solution of the nonlinear problem \eqref{Phit-NSMPB-weak}. When $\|p\|_{H^{1}(\Omega)}$ is sufficiently small, we can use \eqref{b-expansion} to approximate the nonlinear equation  $b(\Phit^{(0)}+p,v)=0$ as a linear equation of $p$,
\[ b(\Phit^{(0)},v) + b^{\prime}(\Phit^{(0)},v;p) = 0,\] 
or equivalently,
\[ b^{\prime}(\Phit^{(0)},v;p) = -b(\Phit^{(0)},v),
 \]
from which we get a linear variational problem for computing $p_0$ as follows: 
\[ \mbox{Find $p_{0} \in {\cal M}_{0}$ such that} \qquad b^{\prime}(\Phit^{(0)},v;p_0) = -b(\Phit^{(0)},v)  \qquad\forall v \in  {\cal M}_{0}. \]
We then can construct a update, $\Phit^{(1)}$, of $\Phit^{(0)}$ by
\[  \Phit^{(1)} = \Phit^{(0)} + p_0.   \]
We next can use mathematical induction to construct a sequence of Newton iterates by the recursive formula
\begin{equation}
\label{newton-phit}
    \Phit^{(k+1)} = \Phit^{(k)} + p_k, \quad k=0,1,2,\ldots,
\end{equation}
where $p_k$ is a solution of the linear variational problem: Find $p_{k} \in {\cal M}_{0}$ such that
\begin{equation}\label{newton-p}
     b^{\prime}(\Phit^{(k)},v;p_{k})=- b(\Phit^{(k)},v) \qquad \forall v\in {\cal M}_{0}.
\end{equation}

In the above construction, we have assumed that the norm $\|p_k\|_{H^{1}(\Omega)}$ is small enough to ensure the convergence of the Newton iterative method. However, in practice, this assumption may not be satisfied, causing a divergence issue. To deal with such an issue, we introduce a damping parameter, $\omega_{k}\in (0, 1]$, to modify the Newton iterative scheme \eqref{newton-phit} into a modified Newton iterative method as follows: 
\begin{equation}
\label{modified_newton-phit}
    \Phit^{(k+1)} = \Phit^{(k)} + \omega_{k} p_k, \quad k=0,1,2,\ldots, 
\end{equation}
where we select $\omega_k$ using a simple scheme as given in Algorithm~1.

{\bf Algorithm~1} (Our damping parameter selection scheme) {\em Let $\eta$ be the smallest damping parameter $\omega_k$ allowable for all $k\geq 0$. In each iteration, we start with $\omega_{k}=1$ and reduce it by half if $\omega_k \geq \eta$ and the following inequality fails,
\begin{equation}
\label{newton-rule}
     \|F(\Phit^{(k+1)})\| \leq  \|F(\Phit^{(k)})\|.
\end{equation}
When $\omega_k < \eta$, we restart the iteration by selecting another initial iterate of $\Phit^{(0)}$. By default, we set $\eta=0.01$. 
}

In Algorithm~1, $\| \cdot \|$ is set as the Euclidean vector norm, and $F$ is a vector function, with which the nonlinear finite element variational problem \eqref{Phit-NSMPB-weak} has been expressed as a system of nonlinear algebraic equations in the vector form
\begin{equation}
\label{Phit_system}
          F(\Phit)=\mathbf{0} \quad \mbox{with } F=(b(\Phit, \varphi_{1}),b(\Phit, \varphi_{2}),\ldots,b(\Phit, \varphi_{N_{h}})),
\end{equation}
where $\varphi_{i}$ denotes the $i$th basis function of ${\cal M}_{0}$, and $N_{h}$ is the total number of interior mesh vertices.

Clearly, $\|F(\Phit^{(k)})\|$ gives an absolute residual error of the $k$th iterate $\Phit^{(k)}$. 
Hence, we say that $\Phit^{(k+1)}$ is a better approximation to a solution of \eqref{Phit-NSMPB-weak} than $\Phit^{(k)}$ if the inequality \eqref{newton-rule} holds; otherwise, we reduce the damping parameter $\omega_k$ to improve the approximation. 

In implementation, we terminate the modified Newton iteration \eqref{modified_newton-phit} if the following termination rule is satisfied:
\begin{equation}
\label{conv-stop}
   \| F(\Phit^{(k)}) \| < \epsilon_r \| F(\Phit^{(0)}) \| + \epsilon_a, 
\end{equation}
where $\epsilon_r$ and $\epsilon_a$ denote the relative and absolute error tolerances, respectively. By default, $\epsilon_r=10^{-8}$ and $\epsilon_a=10^{-8}$. 

\subsection{Reformulation of the modified Newton iterative scheme to avoid direct convolution calculations }

However, solving the linear variational problem \eqref{newton-p} numerically is still expensive, since \eqref{newton-p} involves two convolution terms, $p_{k}*Q_{\lambda}$ and $\Phit^{(k)}*Q_{\lambda}$, whose direct calculation is very costly. To avoid any direct convolution calculations, we develop mathematical techniques and justify them as shown in Theorems~\ref{qk_thm}, \ref{pk_thm}, and \ref{pk_qk_thm}.  
\begin{theorem}
     \label{qk_thm}
     Let $q$ denote the convolution $p*Q_{\lambda}$ with $p \in {\cal M}_{0}$ and $Q_{\lambda}$ being defined in \eqref{Q-def}. Then $q$ is a solution of the linear finite element variational problem: Find $q\in {\cal M}_{0}$ such that 
 \begin{equation}\label{qk-eq2}
   \lambda^2\int_{\Omega}\nabla q(\rr)\cdot\nabla v(\rr)\mathrm{d}\rr + 
     \int_{\Omega}[q(\rr)-p (\rr)]v(\rr)\mathrm{d}\rr  = 0  \qquad \forall v\in {\cal M}_{0}.
 \end{equation}
\end{theorem}
 
 { \em Proof.}
 It has been known that $Q_{\lambda}$ is a Yukawa-type kernel \cite{PhysRevLett.93.108104,ref:hydroelecton} satisfying the equation
 \begin{equation*}
 \label{H-eq-1}
      -\lambda^2\Delta Q_{\lambda}(\rr) + Q_{\lambda}(\rr) =\delta, \quad  \rr \in \mathbb{R}^{3}.
 \end{equation*}
 Doing the convolution of $p$ on both sides of the above equation and using the  multiplication property of 
 convolution, we can get an equation of $q$ with $q=p*Q_{\lambda}$ as follows:
 \begin{equation*}
 \label{phi-def1}
  - \lambda^{2}\Delta q(\rr) +q(\rr) =  p(\rr),\quad  \; \rr \in R^{3}.
 \end{equation*}
 We can then show that $q \in {\cal M}_{0}$ and reformulate the above equation into \eqref{qk-eq2}. This completes the proof.
 
 \begin{theorem}
     \label{pk_thm}
 Let $\zeta^{(k)}=\Phit^{(k)}*Q_{\lambda}$ and $q_k = p_k*Q_{\lambda}$ for $\Phit^{(k)} \in {\cal M}_{0}$ and $p_k \in {\cal M}_{0}$. Then, a sequence of iterates, $\zeta^{(k)}$ for $k=0, 1,2,\ldots$, can be generated by the recursive formula
 \begin{equation}
 \label{newton-uk}
     \zeta^{(k+1)} = \zeta^{(k)}+ \omega_k q_k, \qquad k=0,1,2,\ldots,
 \end{equation}
 where $q_k$ is given in Theorem~\ref{qk_thm}  and  $\zeta^{(0)}$ is an initial guess.  
 \end{theorem}

 { \em Proof.}
 We can obtain the recursive formula \eqref{newton-uk} by doing the convolution with $Q_{\lambda}$ on both sides of \eqref{modified_newton-phit}.

Using the above two theorems, we can reformulate the linear variational problem \eqref{newton-p} into a system of two linear variational problems for computing $p_k$ and $q_k,$ which we present in Theorem~\ref{pk_qk_thm}.

\begin{theorem}
\label{pk_qk_thm}
Let $q_{k} =  p_{k}*Q_{\lambda}$, $\utilde{p}_{k} = (p_{k},q_{k})$, and $\utilde{v} = (v_1, v_2)$. The linear variational problem \eqref{newton-p} can be reformulated as a linear variational system in vector form: Find $\utilde{p}_{k} \in {\cal M}_{0}\times {\cal M}_{0}$ such that
\begin{equation}
 \label{Newton-pq}
    A(\utilde{p}_{k}, \utilde{v}) = L(\utilde{v}) \quad\forall\utilde{v} \in {\cal M}_{0}\times {\cal M}_{0},
 \end{equation}
where $L(\utilde{v})= - l(v_1)$ with the linear form $l(\cdot)$ being defined by 
\begin{eqnarray}
\label{b-def2}
 l(v_1)  &= &\ep \int_{D_{p}}\nabla \Phit^{(k)} (\rr)\cdot\nabla v_1(\rr)\mathrm{d}\rr +\ei \int_{D_{s}}\nabla \Phit^{(k)} (\rr)\cdot\nabla v_1(\rr)\mathrm{d}\rr  \nonumber \\
& &+  (\epsilon_s-\epsilon_\infty)\int_{D_s}\nabla \zeta^{(k)}(\rr) \cdot\nabla v_1(\rr) \mathrm{d}\rr	\\		
& & - \beta \int_{D_{s}}  \frac{ \sum\limits_{i=1}^{n}Z_{i} c_{i}^{b}w_i(\rr) e^{-Z_{i}\Phit^{(k)}(\rr)}}
{1+ \gamma  \frac{\bar{v}^2}{v_0}  \sum\limits_{i=1}^n c_{i}^{b}w_i(\rr)e^{-Z_{i}\Phit^{(k)}(\rr)}} v_1 \mathrm{d}\rr,
 \nonumber
\end{eqnarray}
and $A(\utilde{p}_{k}, \utilde{v})$ is defined by 
\begin{eqnarray}
\label{Newton-A}
 A(\utilde{p}_{k}, \utilde{v}) 
    &= &\ep \int_{D_{p}}\nabla p_{k}(\rr)\cdot\nabla v_1(\rr)\mathrm{d}\rr +\ei \int_{D_{s}}\nabla p_{k}(\rr)\cdot\nabla v_1(\rr)\mathrm{d}\rr \nonumber \\
 & & + \; (\epsilon_s-\epsilon_\infty)\int_{D_s}\nabla q_{k}(\rr)\cdot\nabla v_1(\rr)\mathrm{d}\rr
     +\lambda^2\int_{\Omega}\nabla q_{k}(\rr)\cdot\nabla v_2(\rr)\mathrm{d}\rr, \nonumber \\
   &  &+  \int_{\Omega}[q_{k}(\rr)-p_{k}(\rr)]v_2(\rr)\mathrm{d}\rr   \\
   & & +  \beta \int_{D_s}  \frac{ A_1(\Phit^{(k)}) A_3(\Phit^{(k)}) -   \gamma  \frac{\bar{v}^2}{v_0}  \big[A_2(\Phit^{(k)})\big]^2}{ \big[A_1(\Phit^{(k)})\big]^2} p_k v_1 d\rr. \nonumber
\end{eqnarray}
\end{theorem}

{ \em Proof.}
Using the notation $q_k$ and $\zeta^{(k)}$ and setting $\Phit=\Phit^{(k)}$, $p=p_k$, and $v=v_1$ in \eqref{b-def} and \eqref{b-derivative}, we can get the linear form $l(v_1)$ of \eqref{b-def2}. We can then formulate equation \eqref{newton-p} in linear variational form: Find  $p_{k} \in {\cal M}_{0}$ such that 
\begin{equation}\label{newton-p2}
     a(p_{k},v_1; q_k)=-l(v_1) \qquad \forall v_1\in {\cal M}_{0},
\end{equation}
where $a(p_{k},v_1; q_k)$ is a bilinear form in the expression
\begin{equation}
\label{b-derivative2}
\begin{aligned}
a(p_k,v_1; q_k)  = & \ep \int_{D_{p}}\nabla p_k(\rr)\cdot\nabla v_1(\rr)\mathrm{d}\rr +\ei \int_{D_{s}}\nabla p_k(\rr)\cdot\nabla v_1(\rr)\mathrm{d}\rr   \\
& +  (\epsilon_s-\epsilon_\infty)\int_{D_s}\nabla q_k(\rr) \cdot\nabla v_1(\rr) \mathrm{d}\rr		\\		    
& +  \beta \int_{D_s}  \frac{ A_1(\Phit^{(k)}) A_3(\Phit^{(k)}) -   \gamma  \frac{\bar{v}^2}{v_0}  \big[A_2(\Phit^{(k)})\big]^2}{ \big[A_1(\Phit^{(k)})\big]^2} p_k(\rr) v_1(\rr) d\rr. 
\end{aligned}
\end{equation}
We next set $v=v_2$ to turn \eqref{qk-eq2} as an equation of $p_k$ and $q_k$ as follows: 
\begin{equation}\label{qk-eq3}
   \lambda^2\int_{\Omega}\nabla q_{k}(\rr)\cdot\nabla v_2(\rr)\mathrm{d}\rr + 
     \int_{\Omega}[q_{k}(\rr)-p_{k} (\rr)]v_2(\rr)\mathrm{d}\rr  = 0  \qquad \forall v_2\in {\cal M}_{0}.
\end{equation}
Adding \eqref{newton-p2} with \eqref{qk-eq3}, we can obtain the system \eqref{Newton-pq} for computing $p_k$ and $q_k$. This completes the proof.

With Theorems~\ref{qk_thm}, \ref{pk_thm}, and~\ref{pk_qk_thm}, we now can avoid any direct convolution calculations in the implementation of the modified Newton iterative method \eqref{modified_newton-phit} by the following recursive formulas 
\begin{equation}
\label{Newton_scheme2}
  \Phit^{(k+1)} = \Phit^{(k)} + \omega_k p_k, \quad \zeta^{(k+1)} = \zeta^{(k)}+  \omega_k q_k, \quad k=0,1,2,\ldots,
\end{equation}
where $(p_k, q_k)$ is a numerical solution of the linear variational system \eqref{Newton-pq}, and $\omega_k$ is selected by the scheme defined in Algorithm~1.

A ``blow-up problem'' may occur during a numerical solution of the linear variational system \eqref{Newton-pq}, since both the bilinear form $A(\cdot,\cdot)$ and the linear form $l(\cdot)$ contain exponential functions
\[ e^{-Z_{i}(G+\Psi +\Phit^{(k)})}, \quad i=1,2,\ldots, n,\] which can become a blow-up when a value of $-Z_{i}(G+\Psi +\Phit^{(k)})$ is too large. To avoid this problem, we modify them as $e^{\tau}$ whenever $-Z_{i}(G+\Psi +\Phit^{(k)}) >\tau$. A default value of $\tau$ is 40 in our implementation.

\subsection{A linear NSMPB model}
When $|Z_iu| < 1$, we can use  the Taylor expansion of $ e^{-Z_iu}$ to get
\[  e^{-Z_iu} \approx 1 - Z_i u, \quad i=1,2,\ldots, n.\]
Together with the electro-neutrality condition $ \sum_{i=1}^n Z_ic_i^b =0,$
we can then linearize the nonlinear term of the NSMPB model \eqref{NSMPB-def} by
\[  \frac{ \sum\limits_{i=1}^{n}Z_{i} c_{i}^{b}e^{-Z_{i}u(\rr)}}{1+ \gamma  \frac{\bar{v}^2}{v_0} 
 \sum\limits_{i=1}^n c_{i}^{b}e^{-Z_{i}u(\rr)}}
 \approx
     - \frac{ \sum\limits_{i=1}^{n} Z_{i}^2c_{i}^{b}}{1+ \gamma  \frac{\bar{v}^2}{v_0}   \sum\limits_{j=1}^n c_{j}^{b}} u.
\]
Applying the above linear approximation to \eqref{NSMPB-def}, we can derive a linear NSMPB model as follows:
\begin{equation}
\label{NSMPB-def_linear}
 \left\{\begin{array}{cl}
     -\epsilon_p\Delta u(\rr) =\alpha \sum\limits_{j=1}^{n_{p}}z_{j}  \delta_{\rr_{j}}, & \rr \in D_p,\\
\displaystyle -\epsilon_{\infty}\Delta u(\rr) + \frac{\epsilon_s-\epsilon_\infty}{\lambda^2}[u(\rr)-(u\ast Q_{\lambda})(\rr)]
+ \Upsilon  u(\rr) =0, & \rr \in D_s,\\
\displaystyle u(\s^{-}) = u(\s^{+}), \quad
\epsilon_p\frac{\partial u(\s^{-})}{\partial \nn(\s)} = \epsilon_\infty\frac{\partial u(\s^{+})}{\partial \nn(\s)}+(\epsilon_s-\epsilon_\infty)\frac{\partial (u\ast Q_{\lambda})(\s)}{\partial \nn(\s)}, & \s\in\Gamma,\\
   u (\s) =  g(\s), & \s \in \partial \Omega,
\end{array}
\right.
\end{equation}  
where $\Upsilon$ is a physical constant given by
\[ \Upsilon =  \frac{\kappa^2}{1+ \gamma  \frac{\bar{v}^2}{v_0} \sum\limits_{j=1}^n c_{j}^{b}}.\]
Here $\kappa^2$ is the Debye screening parameter given by
$\kappa^2 = 2 \beta I_s$ with $I_s = \frac{1}{2}  \sum\limits_{i=1}^{n} Z_{i}^2c_{i}^{b}$, which is the ionic strength. 

To avoid the solution singularity, we can construct a solution, $u$, of the linear NSMPB model \eqref{NSMPB-def_linear} by 
\[ u=G+\Psi+\Phit_l, \]
where $G$ and $\Psi$ are given in \eqref{Def_G} and \eqref{Psi_b_nonlocal}, respectively, and $\Phit_l$ denotes a solution of the nonlocal linear interface boundary value problem
\begin{equation}
\label{linear-Phit-uniform-new}
\left\{\begin{array}{cl}
\Delta \Phit_l (\rr) = 0, & \rr\in D_p,\\
-\epsilon_\infty\Delta \Phit_l(\rr) + \frac{\epsilon_s-\epsilon_\infty}{\lambda^2}[\Phit_l - (\Phit_l  \ast Q_{\lambda})]  + \Upsilon \Phit_l
= - \Upsilon [\Psi(\rr) + G(\rr)],  & \rr\in D_s,\\
 \Phit_l (\s^{-}) = \Phit_l (\s^{+}), \quad
\epsilon_p\frac{\partial \Phit_l (\s^{-})}{\partial \nn(\s)} - \epsilon_{\infty}\frac{\partial \Phit_l (\s^{+})}{\partial \nn(\s)} = (\epsilon_s-\epsilon_\infty)\frac{\partial ( \Phit_l  \ast Q_{\lambda})(\s)}{\partial \nn(\s)}, & \s\in\Gamma,\\
\Phit (\s) = 0, & \s \in \partial \Omega.
\end{array}\right.
\end{equation}

Let ${\cal M}_{0}$ be the finite element function space defined in \eqref{M0-def}.
We can obtain a finite element variational problem of the linear interface boundary value problem \eqref{linear-Phit-uniform-new} as follows: Find $(\Phit_{l},\zeta_l) \in {\cal M}_{0}\times {\cal M}_{0}$ such that
\begin{equation}
\label{LNMPBE-eq}
\begin{aligned}
  & \ep \int_{D_{p}}\nabla \Phit_l(\rr)\cdot\nabla v_1(\rr)\mathrm{d}\rr +\ei \int_{D_{s}}\nabla \Phit_l(\rr)\cdot\nabla v_1(\rr)\mathrm{d}\rr  \\
& +  (\epsilon_s-\epsilon_\infty)\int_{D_s}\nabla \zeta_l(\rr) \cdot\nabla v_1(\rr) \mathrm{d}\rr  \\
& + \;  \lambda^2\int_{\Omega}\nabla \zeta_l(\rr)\cdot\nabla v_2(\rr)\mathrm{d}\rr + 
     \int_{\Omega}[\zeta_l(\rr)-\Phit_l (\rr)]v_2(\rr)\mathrm{d}\rr
 + \Upsilon \int_{D_{s}}  \Phit_{l}(\rr) v_{1}(\rr)  \mathrm{d}\rr \\
= &  -\Upsilon \int_{D_{s}} [\Psi(\rr) + G(\rr)] v_{1}(\rr)  \mathrm{d}\rr  \qquad \forall (v_{1}, v_{2}) \in {\cal M}_{0}\times {\cal M}_{0},
\end{aligned}
\end{equation}
where we have set $\zeta_l = \Phit_{l} \ast Q_{\lambda}$ and treated it as an unknown function to avoid any direct convolution calculations.  

\subsection{Initial iterate selections}

The convergence of our modified Newton iterative method \eqref{Newton_scheme2} is contingent upon the selection of initial iterates, $\Phit^{(0)}$ and $\zeta^{(0)}$. In this subsection, we present four selections, called Selections 1 to 4 for clarity. They will be used to elucidate the convergence and robustness of our method in the next section. 

Using the strategy reported in \cite[Eq. (55)]{xie-nonlocal2014}, we obtain Selection 1 as follows:
\begin{equation}
\label{NMPBE-initial1}
   \Phit^{(0)}= \Phit_{local} + \Psi_{local}  - \Psi,
   \quad \zeta^{(0)} = \zeta_{1},
\end{equation}
where $\Psi_{local}$ and $\Phi_{local}$ are component functions of the SMPBE finite element solution of the equations, reported in \cite[Eq. (15)]{nuSMPBE2017} and \cite[Eq. (54)]{nuSMPBE2017}, respectively, $\Psi$ is defined in \eqref{Psi_b_nonlocal}, and $\zeta_{1}$ is a solution to the linear finite element variation problem: Find $\zeta_1  \in {\cal M}_0$ such that
\begin{equation*}
\label{u0-weak}
 \lambda^2\int_{\Omega}\nabla \zeta_1 \cdot\nabla v\mathrm{d}\rr + \int_{\Omega}\zeta_1 v \mathrm{d}\rr
 = \int_{\Omega}\Phit^{(0)}(\rr)v(\rr)\mathrm{d}\rr \qquad\forall v\in {\cal M}_{0}.
\end{equation*}
Since the NSMPB model is reduced to the SMPBE model when $\ei=\es$, it includes the SMPBE model as a special case. In this sense, an SMPBE finite element solution can be a reasonable approximation to an NSMPB finite element solution. Thus, Selection~1 can be a good choice of the initial iteration. We can obtain $\Psi_{local}$ and $\Phit_{local}$ using the SMPBE package reported in \cite{nuSMPBE2017}, and $\Psi$ by the linear finite element method reported in \cite[Eq. 33]{xie-nonlocal2014}.

Given that the linear problem \eqref{LNMPBE-eq} provides a reasonable approximation to the nonlinear problem \eqref{Phit-NSMPB-weak}, we employ its solution to derive Selection~2:
\begin{equation}
\label{NMPBE-initial2}
   \Phit^{(0)}=\Phit_{l}, \qquad \zeta^{(0)}=\zeta_l.
\end{equation}
Clearly, Selection~2 is computationally cheaper than Selection~1 and therefore we set it as the default choice to initialize the modified Newton iterative method \eqref{Newton_scheme2}.  

With a given solution, $\varsigma$, we can simply modify the nonlinear problem \eqref{Phit-NSMPB} as the linear problem: 
\begin{equation*}
\left\{\begin{array}{cl}
\Delta \phi (\rr) = 0, & \rr\in D_p,\\
-\epsilon_\infty\Delta \phi(\rr) + \frac{\epsilon_s-\epsilon_\infty}{\lambda^2}[\phi(\rr) - (\phi  \ast Q_{\lambda})(\rr)] =
\beta \frac{ \sum\limits_{i=1}^{n}Z_{i} c_{i}^{b}w_i(\rr)e^{-Z_{i} \varsigma(\rr)}}
{1+ \gamma  \frac{\bar{v}^2}{v_0}  \sum\limits_{i=1}^n c_{i}^{b}w_i(\rr)e^{-Z_{i} \varsigma(\rr)}}, & \rr\in D_s,\\
 \phi (\s^{-}) = \phi (\s^{+}), \quad
\epsilon_p\frac{\partial \phi (\s^{-})}{\partial \nn(\s)} - \epsilon_{\infty}\frac{\partial \phi (\s^{+})}{\partial \nn(\s)} = (\epsilon_s-\epsilon_\infty)\frac{\partial ( \phi  \ast Q_{\lambda})(\s)}{\partial \nn(\s)}, & \s\in\Gamma,\\
\phi (\s) = 0, & \s \in \partial \Omega.
\end{array}\right.
\end{equation*}
We then reformulate the above problem into the finite element variational system:
Find $(\phi,\zeta) \in {\cal M}_{0}\times {\cal M}_{0}$ such that
\begin{equation}
\label{NMPBE-eq_linear}
\begin{aligned}
  & \ep \int_{D_{p}}\nabla \phi(\rr)\cdot\nabla v_1(\rr)\mathrm{d}\rr +\ei \int_{D_{s}}\nabla \phi(\rr)\cdot\nabla v_1(\rr)\mathrm{d}\rr  
+  (\epsilon_s-\epsilon_\infty)\int_{D_s}\nabla \zeta(\rr) \cdot\nabla v_1(\rr) \mathrm{d}\rr  \\
+ & \; 
   \lambda^2\int_{\Omega}\nabla \zeta(\rr)\cdot\nabla v_2(\rr)\mathrm{d}\rr + 
     \int_{\Omega}[\zeta(\rr)-\phi (\rr)]v_2(\rr)\mathrm{d}\rr  \\
= & \; \beta \int_{D_{s}}  \frac{ \sum\limits_{i=1}^{n}Z_{i} c_{i}^{b} w_i(\rr)e^{-Z_{i} \varsigma(\rr)}}
{1+ \gamma  \frac{\bar{v}^2}{v_0}  \sum\limits_{i=1}^n c_{i}^{b}w_i(\rr)e^{-Z_{i} \varsigma(\rr)}} v \mathrm{d}\rr
 \quad\forall (v_{1}, v_{2}) \in {\cal M}_{0}\times {\cal M}_{0}.
\end{aligned}
\end{equation}

Setting $\varsigma$ as a finite element solution of the linear problem \eqref{linear-Phit-uniform-new}, we can get a solution, $(\phi_1,\zeta_1)$, of the system \eqref{NMPBE-eq_linear}. We then obtain Selection~3:
\begin{equation}
\label{NMPBE-initial3}
   \Phit^{(0)}=\phi_1, \qquad \zeta^{(0)}=\zeta_1.
\end{equation}

The nonlocal problem \eqref{linear-Phit-uniform-new} can be reduced to a local problem by setting $\ei = \es$ as follows:
\begin{equation*}
\label{linear-Phit-uniform-local}
\left\{\begin{array}{cl}
\Delta \Phit_l (\rr) = 0, & \rr\in D_p,\\
-\epsilon_s\Delta \Phit_l(\rr)  + \Upsilon \Phit_l(\rr)
= - \Upsilon [\Psi(\rr) + G(\rr)],  & \rr\in D_s,\\
 \Phit_l (\s^{-}) = \Phit_l (\s^{+}), \quad
\epsilon_p\frac{\partial \Phit_l (\s^{-})}{\partial \nn(\s)} = \epsilon_{s}\frac{\partial \Phit_l (\s^{+})}{\partial \nn(\s)}, & \s\in\Gamma,\\
\Phit (\s) = 0, & \s \in \partial \Omega.
\end{array}\right.
\end{equation*}
Setting $\varsigma$ as a finite element solution of the above problem, which is computationally inexpensive, we can get another solution, $(\phi_2,\zeta_2)$, of the system \eqref{NMPBE-eq_linear}. We then obtain Selection 4:
\begin{equation}
\label{NMPBE-initial4}
   \Phit^{(0)}=\phi_2, \qquad \zeta^{(0)}=\zeta_2.
\end{equation}

\section{The NSMPB finite element program package and numerical results}

We have presented our NSMPB finite element iterative method in Section 3. For clarity, we summarize it in Algorithm~2. 

\vspace{3mm}

{\bf Algorithm~2} (NSMPB Finite Element iterative method) {\em A finite element solution, $u$, of the NSMPB model \eqref{NSMPB-def} is calculated in the following five main steps:
\begin{description}
\item[\hspace*{1cm} {\rm\em Step 1.}] Calculate $G$, $\nabla G$, $\hat{G}$, and $\nabla \hat{G}$ by \eqref{Def_G}, \eqref{nable_G-def},
and \eqref{Def-hatG}.
\item[\hspace*{1cm} {\rm\em Step 2.}] Solve the 
 problem \eqref{Psi_b_nonlocal} for $\Psi$ by the finite element method in \cite[Eq. 33]{xie-nonlocal2014}.
\item[\hspace * {1 cm} {\rm\em Step 3.}] Calculate the initial 
iterations $\Phit^{(0)}$ and $\zeta^{(0)}$ by default selection of \eqref{NMPBE-initial2}.
\item[\hspace*{1cm} {\rm\em Step 4.}] Solve the nonlinear problem \eqref{Phit-NSMPB-weak} for $\Phit$ by the modified Newton  method \eqref{Newton_scheme2}.  
\item[\hspace*{1cm} {\rm\em Step 5.}]  Construct $u$  by the solution decomposition: $u = \Phit+\Psi+ G$.
\end{description} 
}

We implemented Algorithm~2 as an NSMPB finite element program package in Python and Fortran. This package is based on the finite element library of the FEniCS project (Version 2019.1.0) \cite{dolfin2012}, our SMPBE program package \cite{nuSMPBE2017}, and our nonlocal modified PBE program package \cite{xie-nonlocal-solver2012}. To ensure modularity and reusability, the NSMPB program package was designed using object-oriented programming techniques, inheriting all data structures and class methods from the SMPBE program package. For example, the mesh generation module from the SMPBE package has become a part of the NSMPB program package to construct an interface-fitted tetrahedral mesh, $\Omega_h$, for the box domain $\Omega$. Also, the NSMPB program package incorporates SMPBE's class methods to download PDB files from the Protein Data Bank (PDB, {\em http://www.rcsb.org/}) using a four-character PDB identification code, generate the corresponding PQR files using the PDB2PQR package \cite{dolinsky2004pdb2pqr}, and compute electrostatic solvation free energies. A PQR file supplements missing data from the original PDB file, such as hydrogen atoms, atomic charge numbers, and atomic radii. These data are essential for generating a triangular molecular surface mesh $\Gamma_h$ that forms the interface $\Gamma$ between the protein region $D_p$ and the solvent region $D_s$. By default, the CHARMM force field is used, and water molecules are removed from the PDB file when generating the PQR file via the PDB2PQR package. 

Additionally, to further improve computational efficiency, we wrote Fortran subroutines to speed up the calculation of the computationally intensive functions used in the NSMPB finite element scheme. These functions include $G$,  $\nabla G$, $\hat{G}$, and $\nabla \hat{G}$ as well as the residual error vectors $F(\Phit^{(k)})$ defined in \eqref{Phit_system}. We used the Fortran-to-Python interface generator \texttt{f2py} ({\em https://numpy.org/doc/stable/f2py/index.html}) to convert these subroutines into Python modules, enabling direct invocation within Python programs.

We conducted numerical experiments on three proteins (PDB IDs: 1CID, 4PTI, and 1C4K) and a mixture of 0.1 mol/L KNO$_3$ (potassium nitrate) and 0.1 mol/L NaCl (table salt) to validate the convergence of the modified Newton iterative method \eqref{Newton_scheme2} and assess the performance of the NSMPB program package. The four ionic species of the mixture, Cl$^-$, NO$_3^-$, K$^+$, and Na$^+$, were indexed from 1 to 4. Their concentration functions were denoted by $c_1$, $c_2$, $c_3$, and $c_4$; their charge numbers were given by $Z_1 = -1$, $Z_2 = -1$, $Z_3 = 1$, and $Z_4 = 1$; and their bulk concentrations by $c^b_i = 0.1$ mol/L for $i = 1$ to 4. The ionic sizes $v_i$ were estimated using their spherical volumes, 
\[ v_i = \frac{4}{3} \pi r_i^3, \quad i=1,2,3,4,\]
where $r_i$ denotes the ionic radius in \AA. Treating each ion as a hydrated sphere, we obtained radii,
\[
 r_1 = 3.32, \quad r_2 = 3.35, \quad r_3 = 3.58, \quad r_4 = 3.31,
\]
from the website {\em https://bionumbers.hms.harvard.edu/bionumber.aspx?\&id=108517}.
The average ionic size $\bar{v}$ was then calculated as $\bar{v} = 163.715 \, \text{\AA}^3$.

In all numerical experiments, the parameters were fixed as $\ep = 2$, $\es = 80$, $\ei = 1.8$, $\lambda = 15$, and the boundary value function $g = 0$. Linear algebraic systems were solved with absolute and relative residual errors of $10^{-8}$ using the preconditioned generalized minimal residual method (GMRES) with incomplete LU preconditioning. For the modified Newton iterative method \eqref{Newton_scheme2}, we set $\tau = 40$, $\eta = 0.01$, $\epsilon_a = 10^{-8}$, and $\epsilon_r = 10^{-8}$. All computations were performed on an Apple M4 Pro chip on our Mac mini computer with 64 GB of memory.

\subsection{Construction of six box meshes per protein case}

\begin{figure}[h!]

         \begin{subfigure}[b]{0.32\textwidth}
                \centering
               \includegraphics[width=0.8\textwidth]{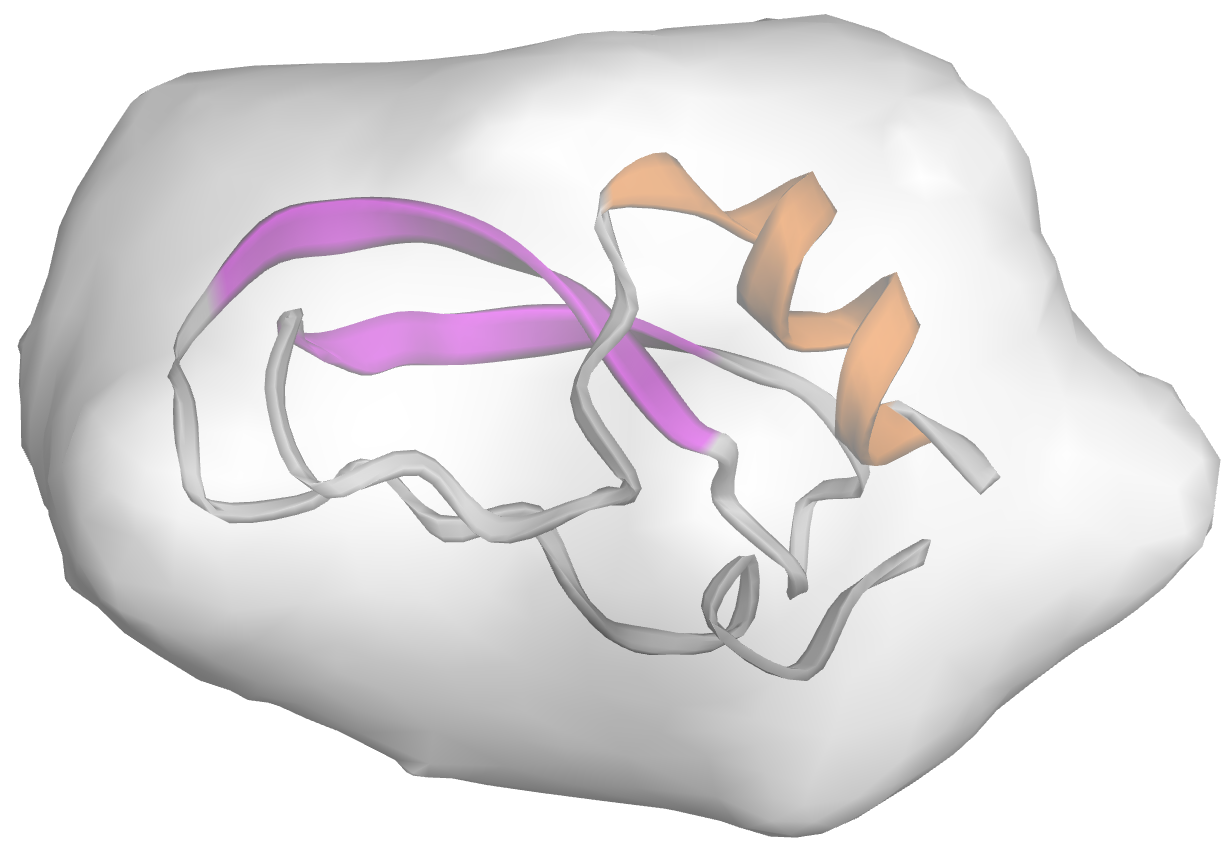}
                \caption{Protein 4PTI in cartoon}
               \label{subfig:4pti_a}
        \end{subfigure}
        \begin{subfigure}[b]{0.32\textwidth}
                \centering
               \includegraphics[width=0.8\textwidth]{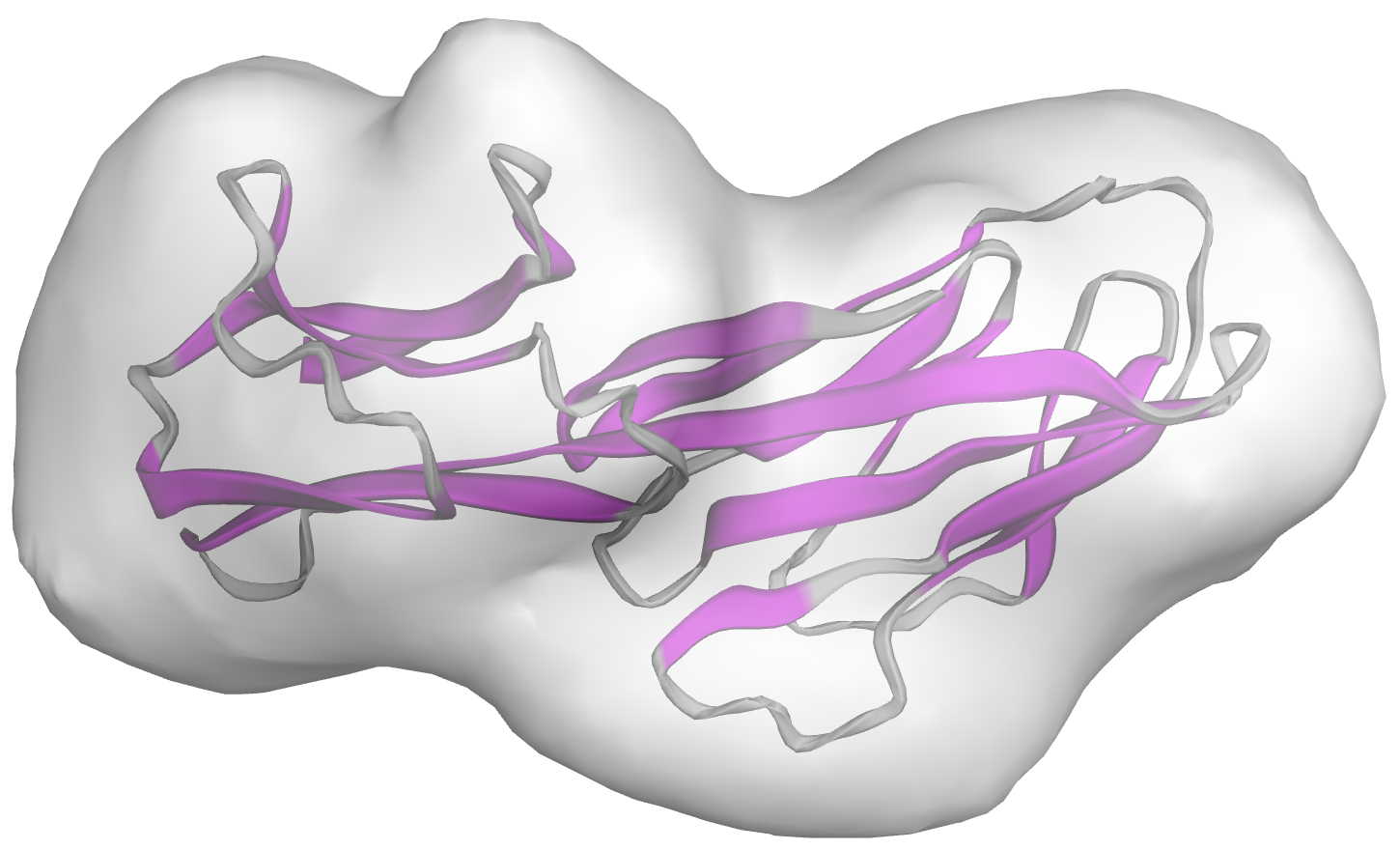}
               \vspace{0.4cm}
                \caption{Protein 1CID in cartoon}
               \label{subfig:1CID_a}
        \end{subfigure} 
        \begin{subfigure}[b]{0.32\textwidth}
                \centering
              \includegraphics[width=0.8\textwidth]{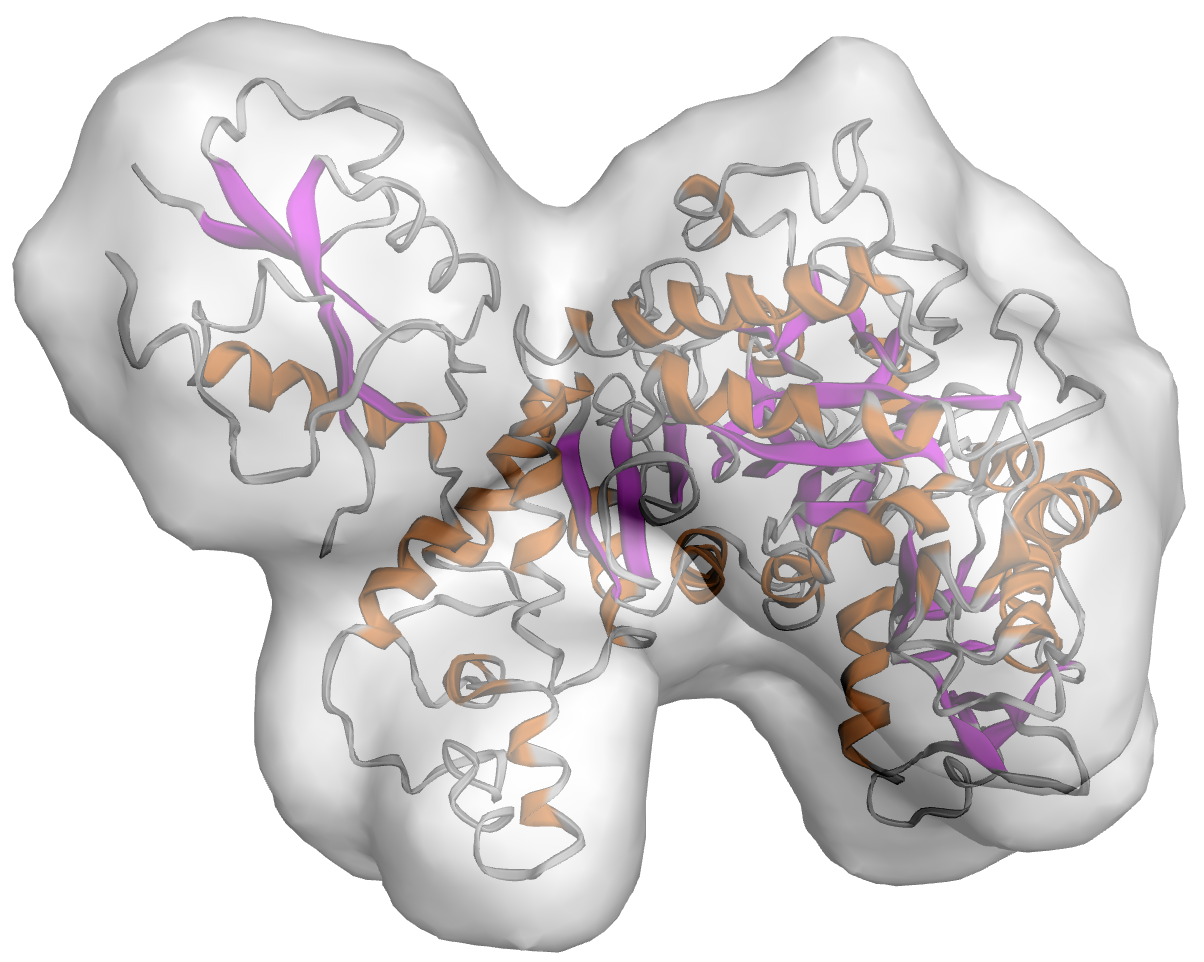}
                \caption{Protein 1C4K in cartoon }
               \label{subfig:solvent-domain-new}
        \end{subfigure}

         \begin{subfigure}[b]{0.32\textwidth}
                \centering
               \includegraphics[width=0.8\textwidth]{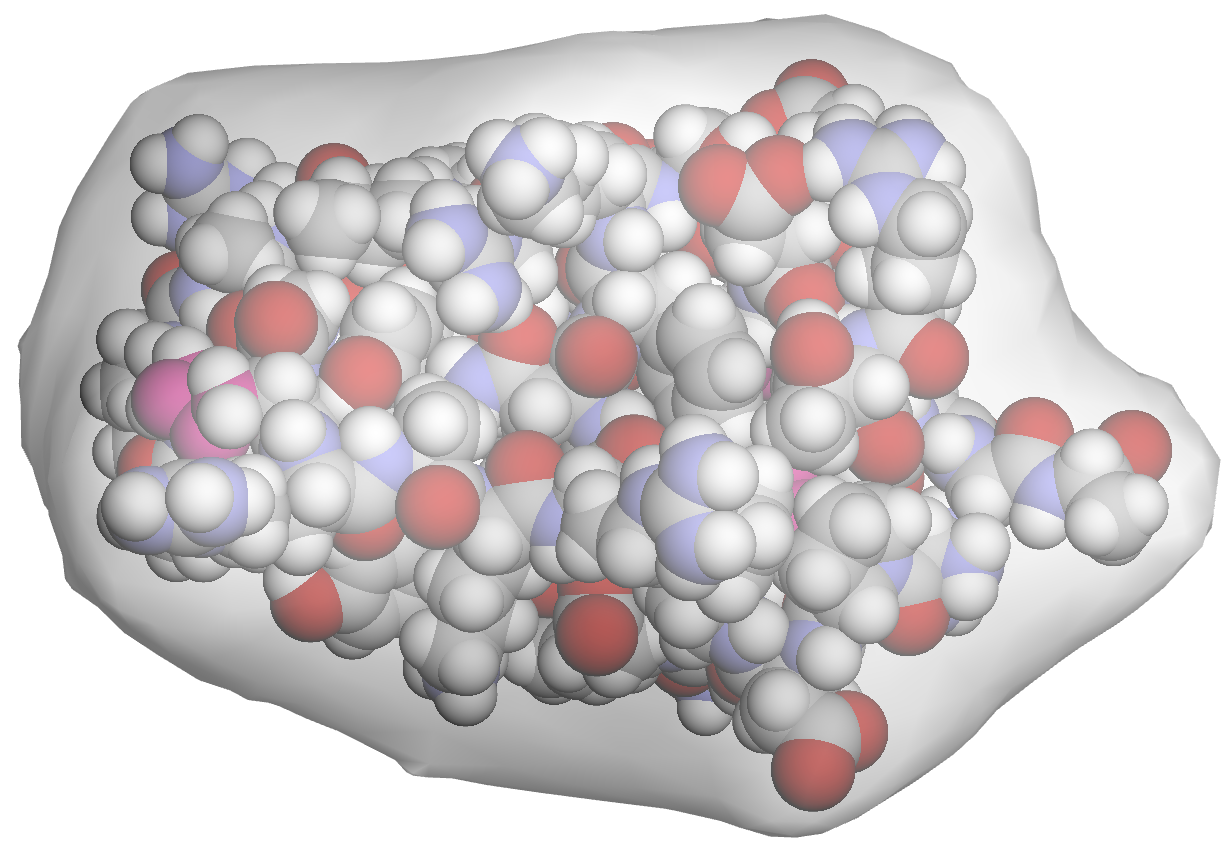}
                \caption{Protein 4PTI in sphere}
               \label{subfig:4pti_b}
        \end{subfigure}
        \begin{subfigure}[b]{0.32\textwidth}
                \centering
               \includegraphics[width=0.8\textwidth]{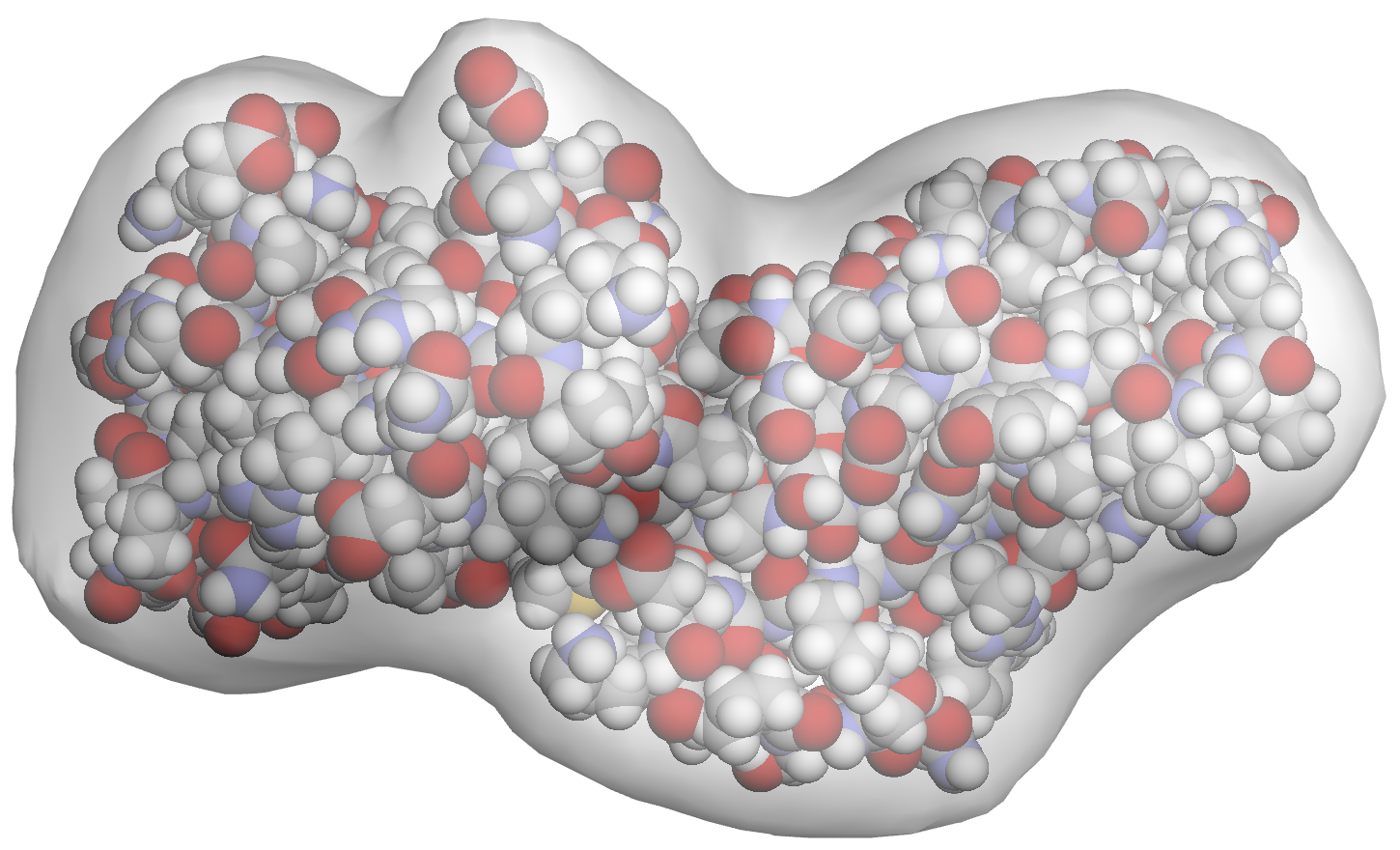}
               \vspace{0.4cm}
                \caption{Protein 1CID in sphere}
               \label{subfig:1CID_b}
        \end{subfigure} 
        \begin{subfigure}[b]{0.32\textwidth}
                \centering
              \includegraphics[width=0.8\textwidth]{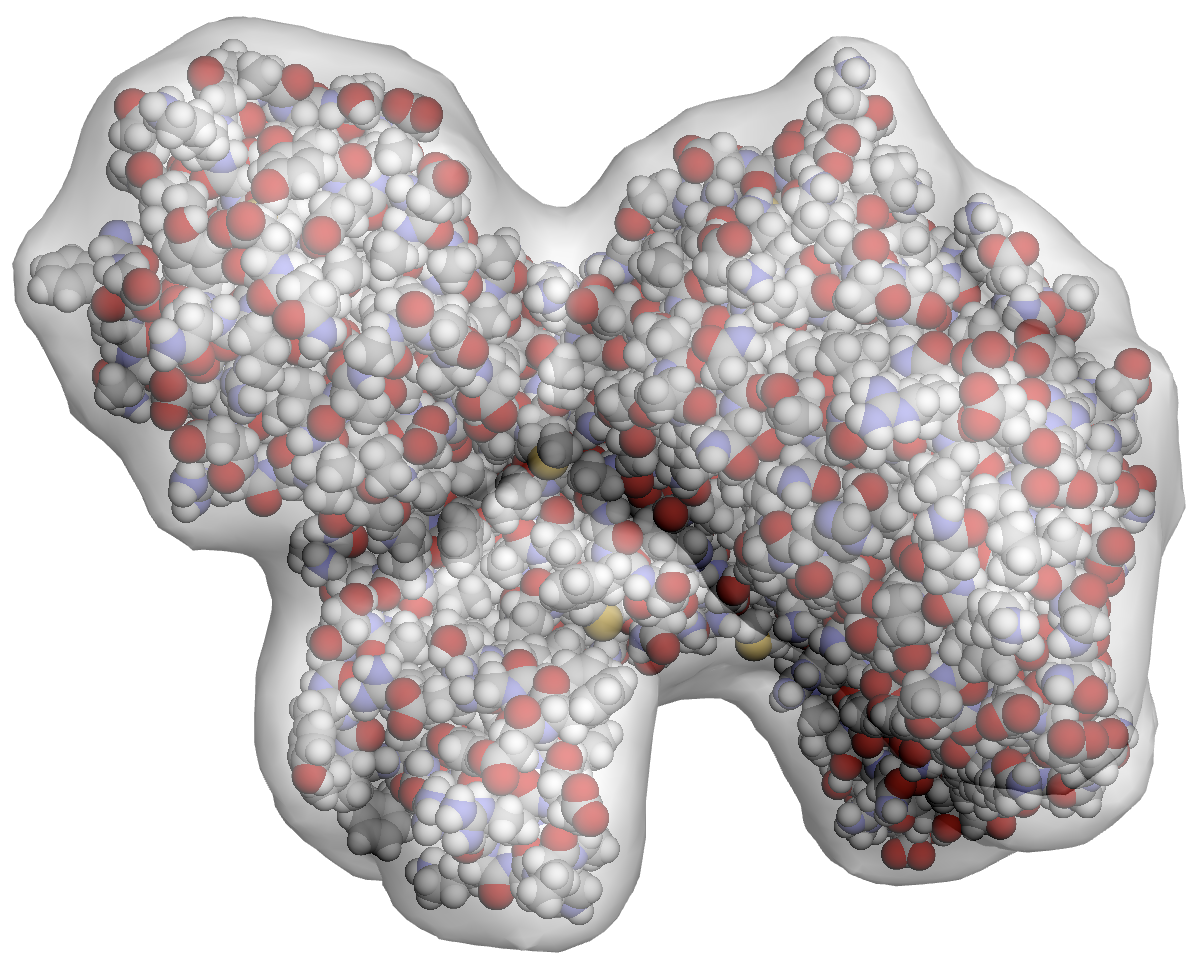}
                \caption{Protein 1C4K in sphere}
               \label{subfig:1C4K_b}
        \end{subfigure}
 \caption{Three proteins' crystallographic molecular structures depicted in cartoon and sphere representations wrapped by the three protein regions  $D_p$ (in gray color) generated by our NSMPB software package. Here 4PTI, 1CID, and 1C4K are the  Protein Data Bank identifications (PDB ID).}       
\label{protein_structures}          
\end{figure} 

With the PDB file of each protein, we determined the smallest box, denoted by $[p_{x1}, p_{x2}] \times [p_{y1}, p_{y2}] \times [p_{z1}, p_{z2}]$, that contains the protein. We then constructed the box domain $\Omega$ by
\[
\Omega = \{ (x, y, z) \in \R^3 \mid L_{x1} < x < L_{x2}, \; L_{y1} < y < L_{y2}, \; L_{z1} < z < L_{z2} \},
\]
where $L_{x1} = p_{x1} - a_1$, $L_{x2} = p_{x2} + a_1$, $L_{y1} = p_{y1} - a_2$, $L_{y2} = p_{y2} + a_2$, $L_{z1} = p_{z1} - a_3$, and $L_{z2} = p_{z2} + a_3$. We set $a_i = 30$ for $i = 1, 2, 3$ for simplicity. Such a box contains the protein in its central part. The dimensions of the three box domains constructed for the three proteins are given in Table~\ref{table: mesh-data}.

The mesh generation module of the NSMPB package incorporates the TMSmesh package (Version 2.1) \cite{chen2012triangulated,liu2018efficient} to generate a molecular triangular surface mesh, $\Gamma_h$. For each protein provided in a PQR file, we generated a molecular triangular surface mesh, $\Gamma_h$, setting the TMSmesh parameters $d$, $c$, and $e$ to 0.1, 0.9, and 0.9, respectively. This mesh served as the interface $\Gamma$ between the protein region $D_p$ and the solvent region $D_s$.

Figure~\ref{protein_structures} presents the three protein regions (in gray) generated by the mesh module of our NSMPB package. It reveals that these protein regions possess irregular shapes and effectively encapsulate the crystallographic molecular structures of these three proteins depicted in cartoon and sphere representations. 

We constructed uniform triangular meshes on the six side surfaces of the box domain $\Omega$ dividing the intervals $[L_{x1}, L_{x2}]$, $[L_{y1}, L_{y2}]$, and $[L_{z1}, L_{z2}]$ into 10 equal subintervals each. Using this box surface meshes and the interface mesh $\Gamma_h$, we generated an irregular interface-fitted tetrahedral mesh, $\Omega_h$, of $\Omega$, called Mesh 1, by the TetGen software (Version 1.5) \cite{si2015tetgen} such that the box mesh  $\Omega_h$ satisfies the partition:
\begin{equation}
\label{box_partition}
 \Omega_h = D_{p,h} \cup D_{s,h} \cup \Gamma_h,
\end{equation}
where $D_{p,h}$ and $D_{s,h}$ are the irregular tetrahedral meshes of $D_p$ and $D_s$, respectively, fitting on $\Gamma_h$.  Such an interface-fitted tetrahedral mesh can significantly improve the numerical accuracy of our NSMPB finite element solver.
 
We then constructed five finer meshes than Mesh~1, named Meshes 2 to 6, using the TetGen parameter {\em -a} equal to 100, 10, 5, 3, and 1, which restrict the maximum tetrahedron volumes of Meshes 2 to 6 not exceeding 100, 10, 5, 3, and 1, respectively. We also used the TetGen mesh quality parameter {\em -q} with the value 1.2 (that is, the maximum radius-edge ratio of the tetrahedra is limited to 1.2) to ensure that these five finer meshes are of high quality. The mesh data of these meshes are reported in Table~\ref{table: mesh-data}. 

\begin{table}[h!]
\centering
\scalebox{1}{
  \begin{tabular}{|c||c|c|c|c|c||c|c|c|}
   \hline
Mesh  & \multicolumn{5}{c||}{Number of vertices}&\multicolumn{3}{c|}{Number of tetrahedra } \\  \cline{2-9}
   &$\Omega_h$&  $D_{s,h}$ &  $D_{p,h}$ & $\Gamma_{h}$ &$\partial\Omega_h$ &$\Omega_h$&  $D_{s,h}$ &  $D_{p,h}$ \\
\hline\hline
    \multicolumn{9}{|c|}{Protein 4PTI with 892 atoms in $\Omega = [-23, 61]\times [-12, 66] \times [-40, 54]$}  \\ \cline{1-9}
\hline\hline
Mesh 1 & 2957 & 2668 & 1588 & 1299 & 602 & 15825 & 10512 & 5313 \\ \hline 
Mesh 2 & 10041 & 8778 & 2571 & 1308 & 2348 & 55104 & 43361 & 11743 \\ \hline 
Mesh 3 & 41971 & 39192 & 5541 & 2762 & 10923 & 230232 & 204515 & 25717 \\\hline 
Mesh 4 & 74074 & 68423 & 10611 & 4960 & 14556 & 423573 & 373034 & 50539 \\ \hline 
Mesh 5 & 116285 & 107231 & 16300 & 7246 & 19813 & 676381 & 597294 & 79087 \\ \hline 
Mesh 6 & 310737 & 291226 & 34316 & 14805 & 43568 & 1843368 & 1674986 & 168382 \\ \hline 
 \hline\hline
  \multicolumn{9}{|c|}{Protein 1CID with 2783 atoms in $\Omega = [-60, 40] \times [-13, 88] \times [-27, 94]$}  \\
    \cline{1-9}
Mesh 1 & 5439 & 4732 & 3450 & 2743 & 602 & 30757 & 18763 & 11994 \\ \hline
Mesh 2 & 17601 & 14565 & 5812 & 2776 & 3247 & 99467 & 72150 & 27317 \\ \hline 
Mesh 3 & 65288 & 58557 & 12730 & 5999 & 12938 & 372037 & 311351 & 60686 \\ \hline 
Mesh 4 & 122568 & 108737 & 24534 & 10703 & 19445 & 715803 & 596086 & 119717 \\ \hline 
Mesh 5 & 193592 & 172063 & 37401 & 15872 & 27007 & 1143720 & 959546 & 184174 \\ \hline
Mesh 6 & 498790 & 452577 & 78166 & 31953 & 55708 & 3003448 & 2613291 & 390157 \\ \hline
   \hline\hline
    \multicolumn{9}{|c|}{Protein 1C4K with 11439 atoms in $\Omega = [4, 154] \times [-42, 97] \times [-32, 108]$}  \\ \cline{1-9}
\hline\hline
Mesh 1 & 14590 & 12215 & 10427 & 8052 & 602 & 85715 & 48201 & 37514 \\ \hline
Mesh 2 & 45068 & 34518 & 18753 & 8203 & 5061 & 266757 & 175759 & 90998 \\ \hline 
Mesh 3 & 150109 & 130059 & 36255 & 16205 & 20920 & 885450 & 709328 & 176122 \\ \hline
Mesh 4 & 296086 & 255047 & 70038 & 28999 & 37537 & 1759132 & 1411379 & 347753 \\ \hline
Mesh 5 & 459491 & 396694 & 104741 & 41944 & 48619 & 2765758 & 2240437 & 525321 \\ \hline
Mesh 6 & 1188840 & 1046984 & 228697 & 86841 & 96280 & 7278073 & 6112150 & 1165923 \\ \hline
\hline
  \end{tabular}} 
  \caption{Mesh data for the three protein cases. Here $\Omega_h$,  $D_{s,h}$,  $D_{p,h}$, $\Gamma_{h}$, and $\partial\Omega_h$ denote the meshes of the box domain $\Omega$,  solvent region $D_{s}$, protein region $D_{p}$,  boundary $\Gamma$ of $D_p$, and boundary $\partial\Omega$ of $\Omega$, respectively. }
   \label{table: mesh-data}
\end{table} 

\begin{figure}[h!]

         \begin{subfigure}[b]{0.32\textwidth}
                \centering
               \includegraphics[width=0.85\textwidth]{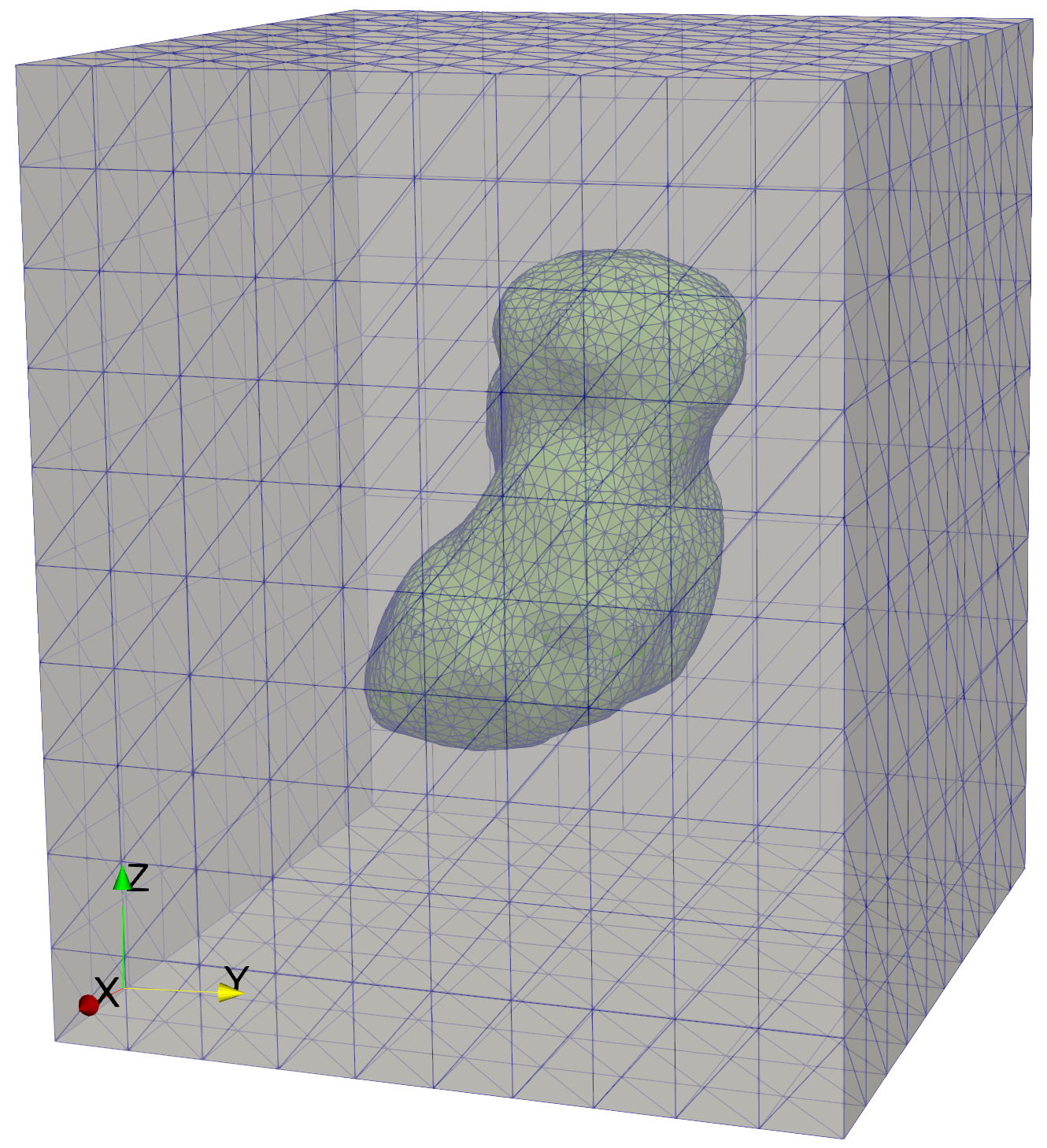}
                \caption{Mesh 1}
               \label{subfig:1C1D_mesh1}
        \end{subfigure}
        \begin{subfigure}[b]{0.32\textwidth}
                \centering
               \includegraphics[width=0.85\textwidth]{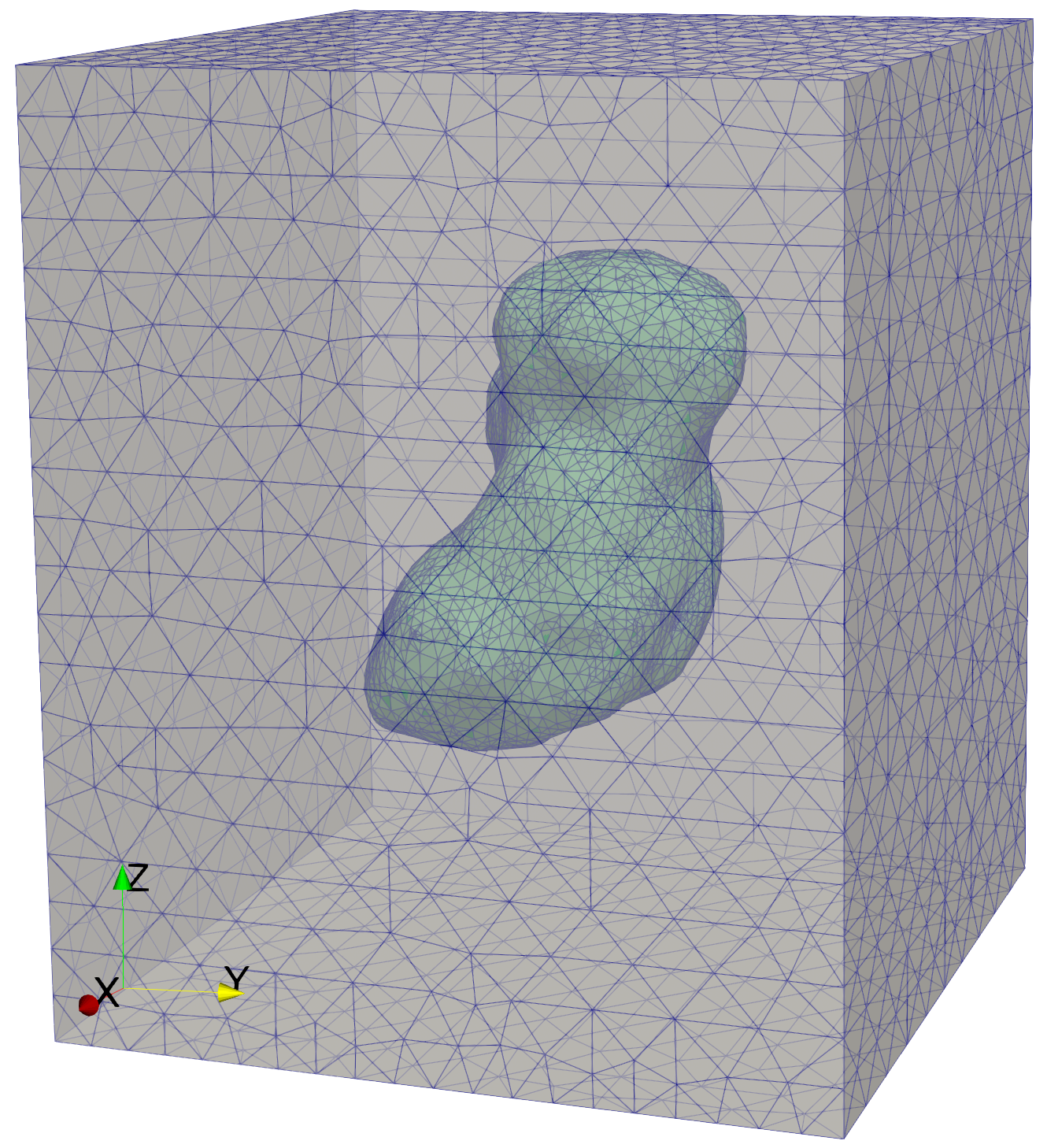}
                \caption{Mesh 2}
               \label{subfig:1C1D_mesh2}
        \end{subfigure} 
        \begin{subfigure}[b]{0.32\textwidth}
                \centering
              \includegraphics[width=0.85\textwidth]{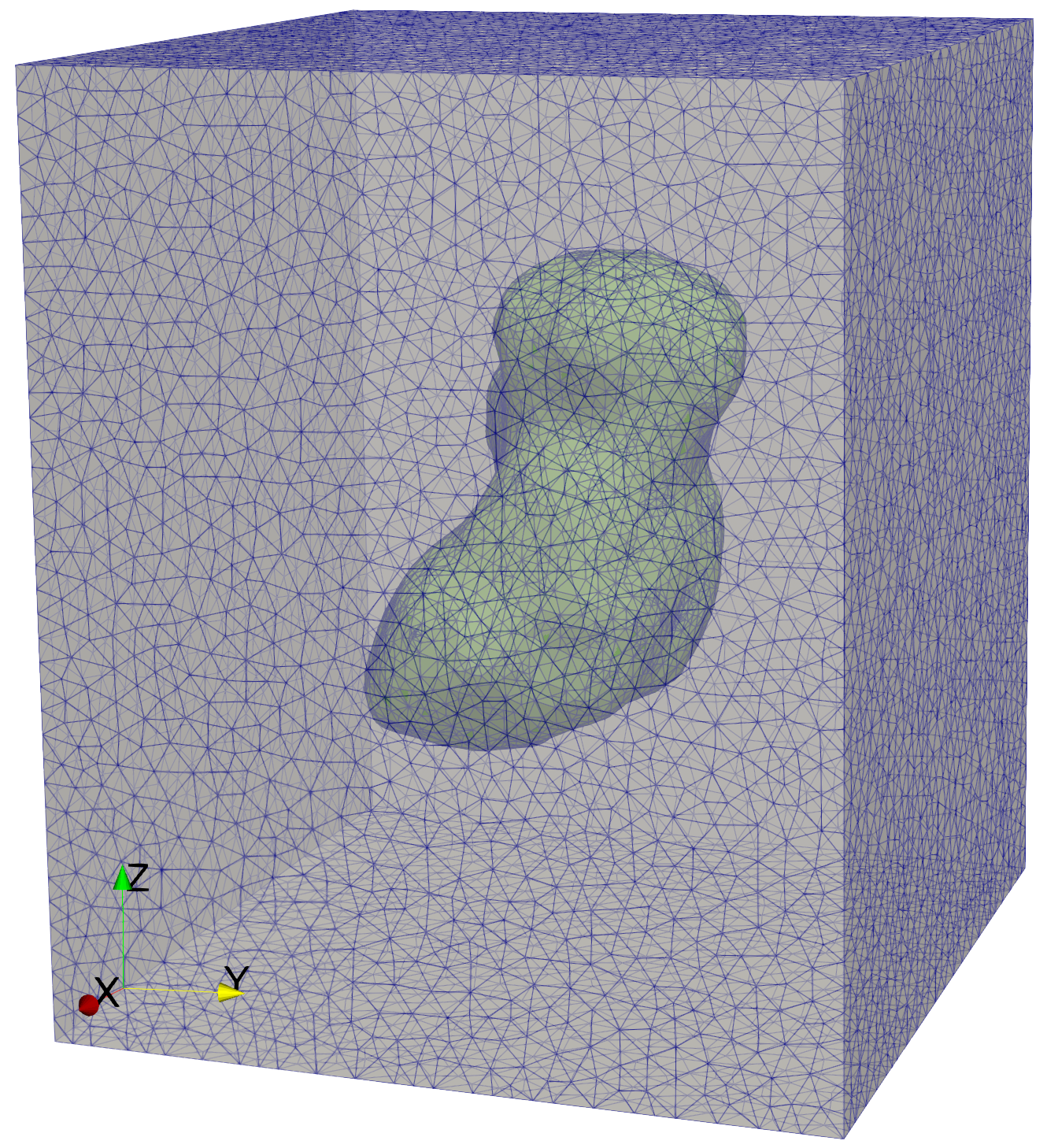}
                \caption{Mesh 3 }
               \label{subfig:1C1D_mesh3}
        \end{subfigure}

\vspace{0.2cm}

         \begin{subfigure}[b]{0.32\textwidth}
                \centering
               \includegraphics[width=0.8\textwidth]{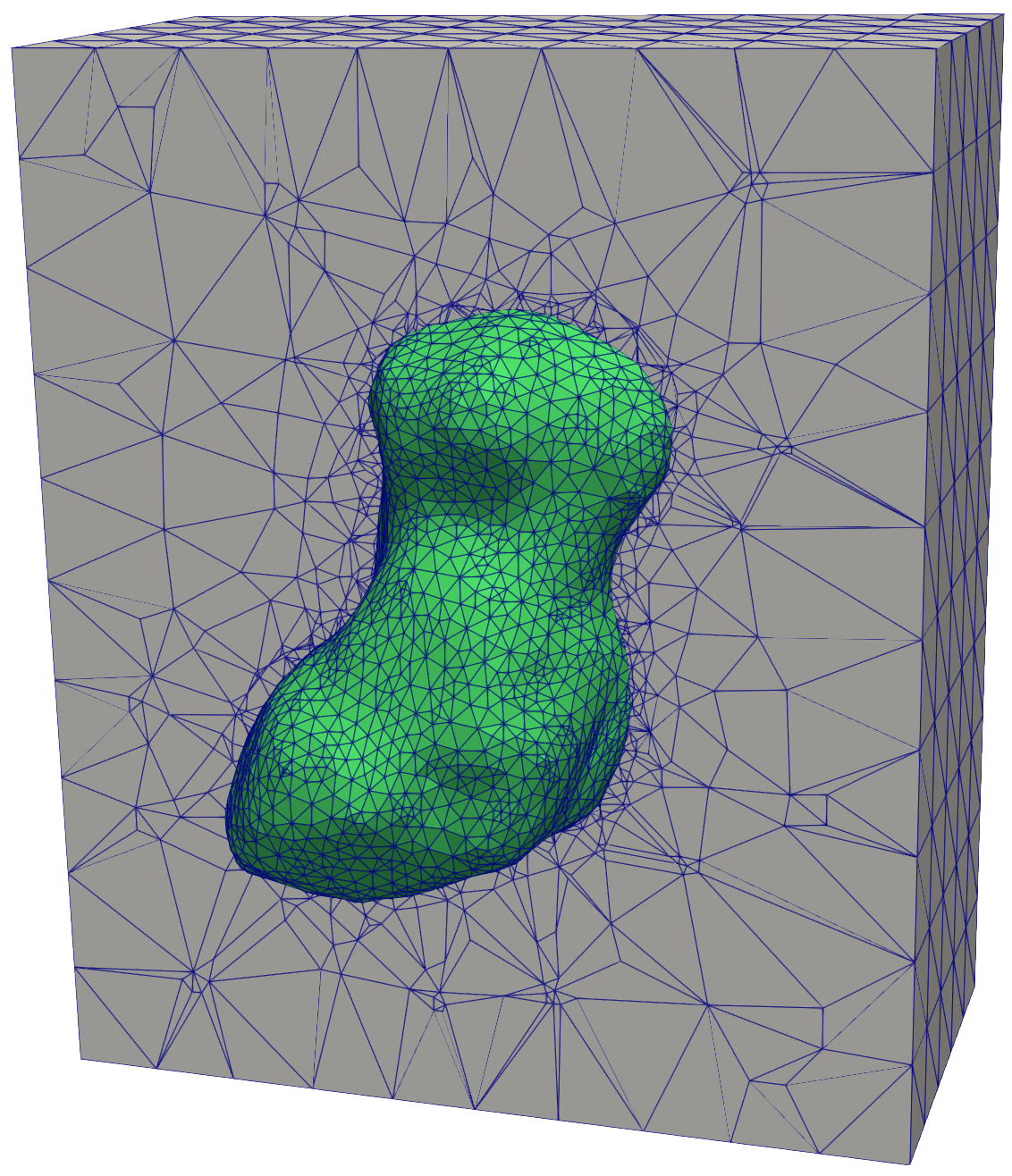}
                \caption{A clip view  of Mesh 1}
               \label{subfig:1C1D_clip1}
        \end{subfigure}
        \begin{subfigure}[b]{0.32\textwidth}
                \centering
               \includegraphics[width=0.8\textwidth]{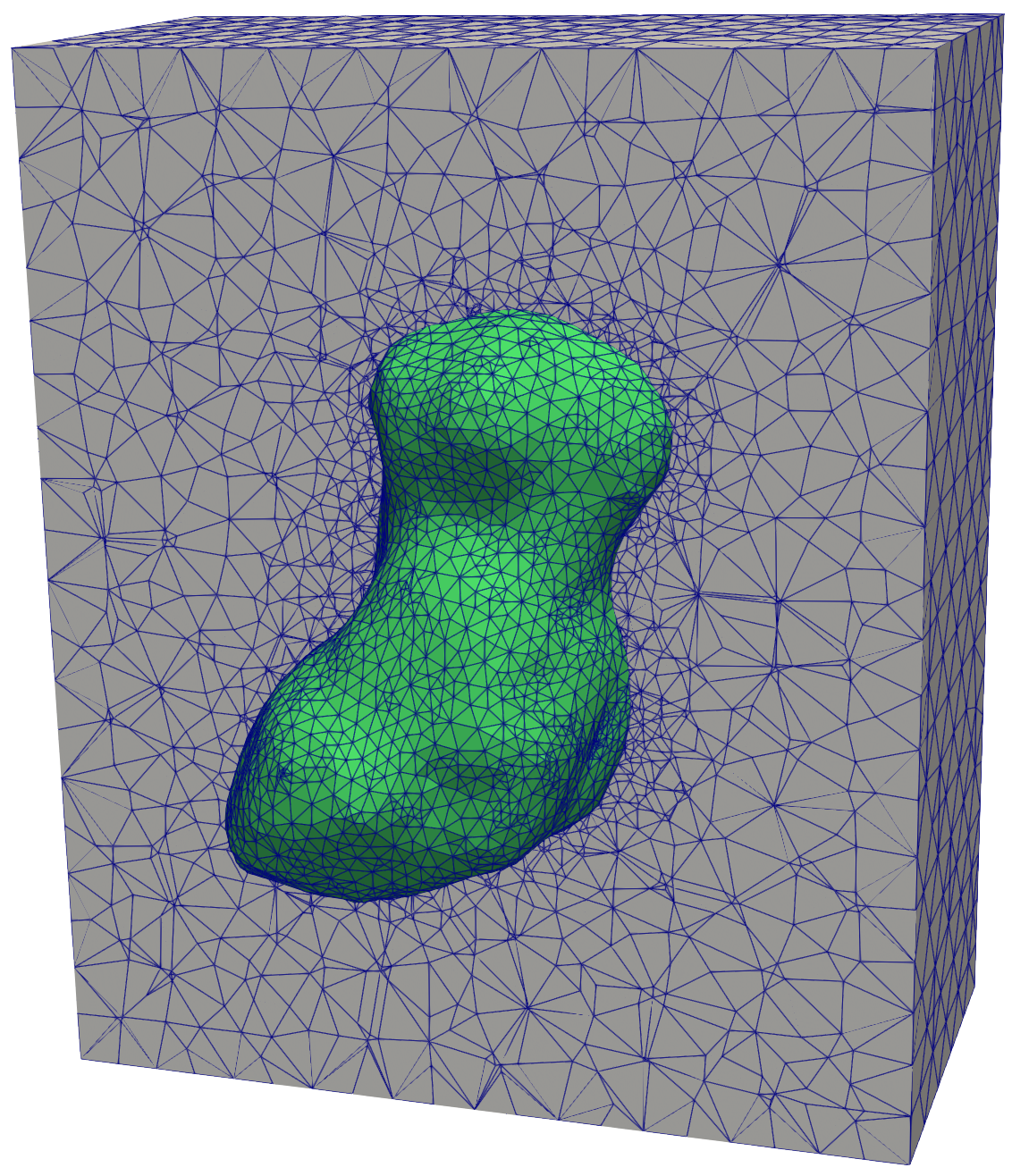}
                \caption{A clip view of Mesh 2}
               \label{subfig:1C1D_clip2}
        \end{subfigure} 
        \begin{subfigure}[b]{0.32\textwidth}
                \centering
              \includegraphics[width=0.8\textwidth]{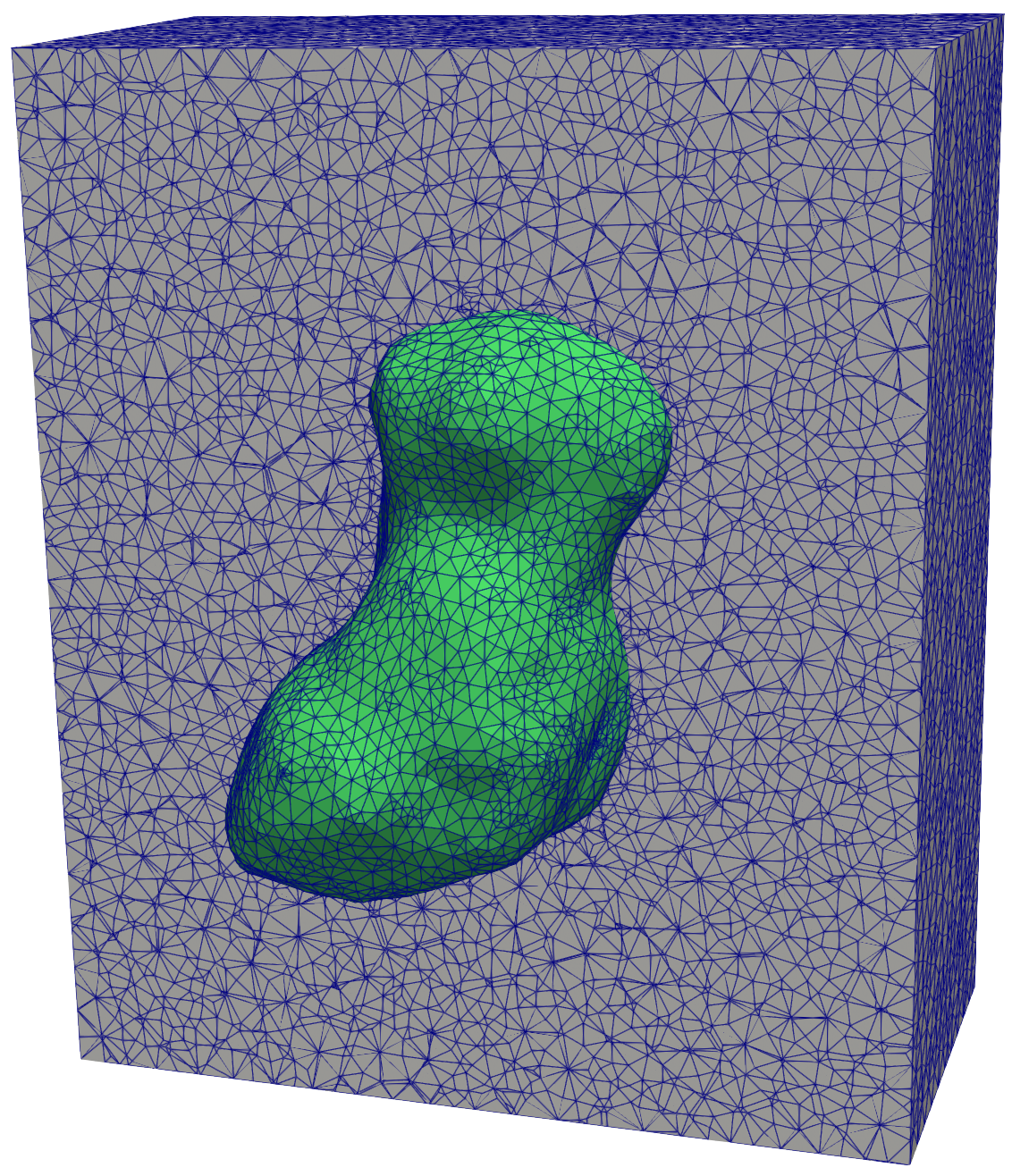}
                \caption{A clip view  of Mesh 3}
               \label{subfig:1C1D_clip3}
        \end{subfigure}
 \caption{(a, b, c) One view of the box mesh $\Omega_{h_j}$ (i.e., Mesh $j$ for $j=1,2,3$) generated by our NSMPB software package in the protein 1CID case. (d, e, f) One clipped view (with $x=0$) of the solvent region mesh $D_{s,h_j}$, along with the triangular boundary surface mesh of the protein region mesh $D_{p, h_j}$, shown in green, highlighting the mesh irregularity between  $D_{s,h_j}$ and  $D_{p, h_j}$. }       
\label{protein_meshes}          
\end{figure} 

Table~\ref{table: mesh-data} shows that the number of mesh vertices and the number of tetrahedra within the protein region $D_p$ and the solvent region $D_s$, as well as the number of mesh vertices on the box boundary $\partial \Omega$, have increased dramatically from Mesh 1 to Mesh 6, due to increasingly strict constraints on the volumes of tetrahedra. The number of vertices on the mesh of the interface surface $\Gamma_h$ does not change significantly from Mesh 1 to 6, because the mesh of the interface surface of Mesh 1 is already quite dense. 

Figure~\ref{protein_meshes} presents two views of Meshes 1 to 3 in the protein 1CID case, as examples, to illustrate mesh refinement processes and mesh irregularities. Here, the upper row of Figure~\ref{protein_meshes} shows the box meshes, $\Omega_{h_j}$, with low opacity, to make the protein region meshes $D_{p, h_j}$ visible. In contrast, the lower row shows clipped views of the solvent region meshes, $D_{s, h_j}$, providing a clearer view of $D_{p, h_j}$ and the irregularities near the interface mesh $\Gamma_{h_j}$. From Figure~\ref{protein_meshes} we can also see that the interface mesh $\Gamma_{h_1}$ of Mesh 1 has been significantly refined with a much higher density compared to the solvent region mesh $D_{s,h_1}$ near the boundary $\partial \Omega$ of the box domain $\Omega$. Thus, in Meshes 2 and 3, relatively little refinement occurs in $\Gamma_{h_2}$ and $\Gamma_{h_3}$, respectively, when the maximum tetrahedral volume constraint parameter ({\em -a}) of TetGen is applied to generate Meshes 2 and 3 --- the two refinements of Mesh 1.  

\begin{table}[h!]
\centering
\begin{tabular}{|c||c|c|c|c|c||c|}
\hline
 Mesh  & Generate & Calculate $G$, & Find $\Psi$ 	& Find $\Phit^{(0)}$& Find $\Phit$ &  Total\\
 & meshes & $\hat{G}, \nabla G, \nabla \hat{G}$	& by \cite{xie-nonlocal2014}	&  by \eqref{NMPBE-initial2}		& by \eqref{Newton_scheme2}	 & CPU time\\
    \hline\hline
    \multicolumn{7}{|c|}{Cases of Protein 1CID with 2783 atoms in $\Omega = [-60, 40] \times [-13, 88] \times [-27, 94]$}  \\
    \cline{1-7}
\hline
Mesh 1  & $1.09 $		& $0.05$		& $0.12$	& $ 0.12$			& $1.63$		 & $5.04$\\
Mesh 2   & $1.57$		& $0.15$		& $0.42$	& $0.37$			& $4.67$	& $10.41$\\
Mesh 3  & $4.45$		& $0.58$		     & $2.15$	    & $1.9$			    & $21.63$	& $39.32$\\
Mesh 4     & $8.35$		& $1.18$		& $5.25$	& $4.36$			& $43.62$ & $1.29$ min.\\
Mesh 5 	& $12.75$		& $1.79$		& $8.66$	& $7.17$		& $1.22$ min.	& $2.08$ min. \\
Mesh 6 	& $32.92$		& $4.52$		& $29.41$	& $23.31$		& $3.87$ min.	& $6.15$ min.\\
\hline\hline
    \multicolumn{7}{|c|}{Cases of Protein 4PTI with 892 atoms in $\Omega = [-23, 51]\times [-13, 56] \times [-30, 44]$}  \\ \cline{1-7}
\hline
Mesh 1  & $0.64 $		& $0.01$		& $0.07$	& $ 0.06$			& $0.85$		 & $2.70$\\
Mesh 2 & $1.15$		& $0.03$		& $0.23$	& $0.2$			& $2.69$	& $6.08$\\
Mesh 3	& $2.92$		& $0.13$		& $1.14$	& $1.08$			& $12.38$	& $22.93$\\
Mesh 4	& $4.72$		& $0.22$		& $2.62$	& $2.32$			& $27.42$	& $45.84$\\
Mesh 5 & $7.30$		& $0.35$		& $4.64$	& $4.1$		& $41.86$	& $1.17$ min.\\
Mesh 6 & $20.40$		& $0.92$		& $16.23$	& $13.58$		& $2.07$ min.	& $3.37$ min.\\
\hline\hline
    \multicolumn{7}{|c|}{Cases of Protein 1C4K with 11439 atoms in $\Omega = [4, 154] \times [-42, 97] \times [-32, 108]$}  \\ \cline{1-7}
\hline
Mesh 1  & $2.33 $		& $0.52$		& $0.33$	& $0.30$			& $5.33$		 & $14.66$\\
Mesh 2   & $3.76$		& $1.61$		& $1.4$	& $1.16$			& $15.46$	& $32.65$\\
Mesh 3	& $10.10$		& $5.43$		& $6.16$	& $5.01$			& $57.57$	& $1.87$ min.\\
Mesh 4  	     & $20.34$		& $11.19$		& $15.49$	& $11.57$			& $2.36$ min.	& $4.21$ min.\\
Mesh 5 	& $31.26$		& $17.44$		& $26.32$	& $19.85$		& $3.63$ min.	& $6.50$ min.\\
Mesh 6 	& $88.80$ 		& $44.28$		& $1.80$ min.	& $72.60$  & $12.37$ min. & $20.15$ min.\\
\hline
\end{tabular}
\caption{Performance of the NSMPB program package for the three proteins and a mixture of 0.1 mole KNO$_3$ and 0.1 mole NaCl in computer CPU times in seconds except where noted. Mesh data are given in Table~1. These tests were done on our Mac mini computer with Apple M4 chip and 64 GB memory.
}
\label{Performance-CPU}
\end{table}

\begin{figure}[h!]
    \centering
    \begin{subfigure}[b]{0.315\textwidth}
                \centering
               \includegraphics[width=\textwidth]{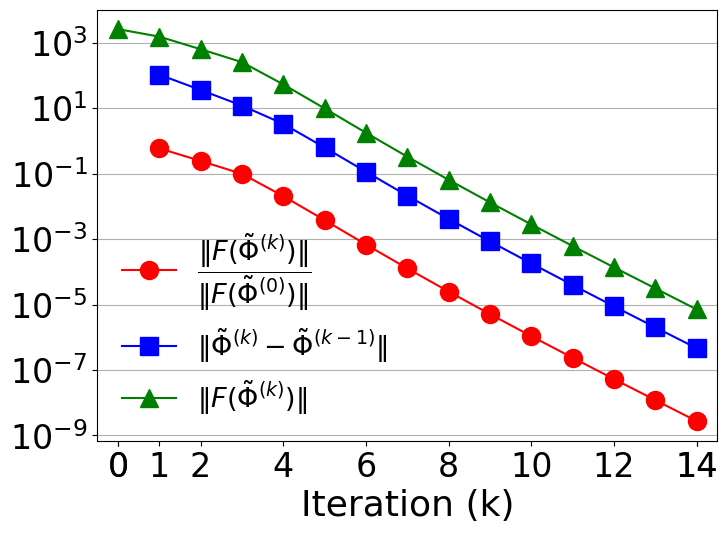}
                \caption{Case of Mesh 1}
               \label{subfig:mesh1-conv}
        \end{subfigure}
        \hspace{0.1cm} 
        \begin{subfigure}[b]{0.315\textwidth}
                \centering
               \includegraphics[width=\textwidth]{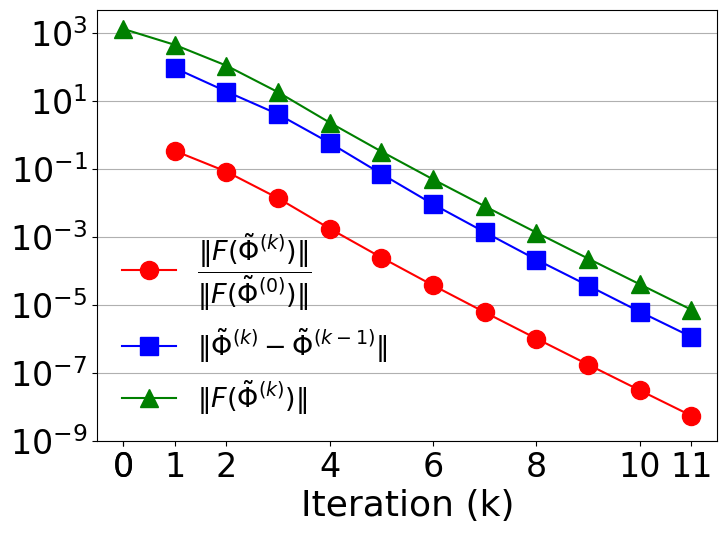}
                \caption{Case of Mesh 2}
               \label{subfig:mesh2-conv}
        \end{subfigure} 
        \hspace{0.1cm}
        \begin{subfigure}[b]{0.315\textwidth}
                \centering
              \includegraphics[width=\textwidth]{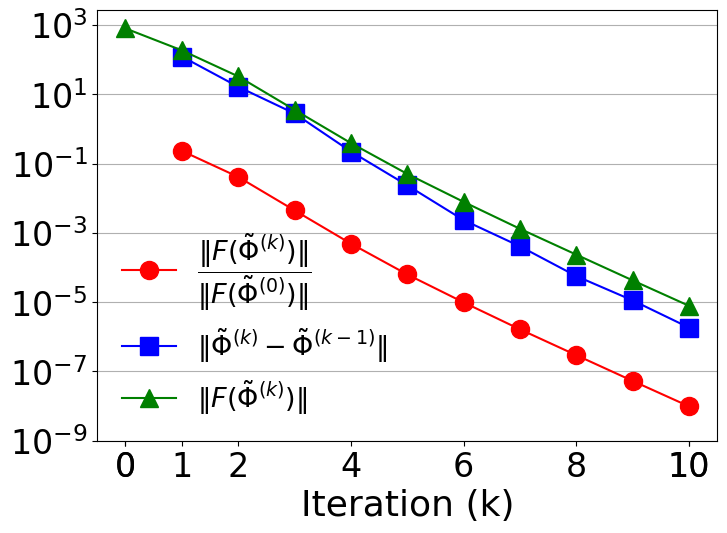}
                \caption{Case of Mesh 3}
               \label{subfig:mesh3-conv}
        \end{subfigure}
        
        
        \begin{subfigure}[b]{0.315\textwidth}
                \centering
              \includegraphics[width=\textwidth]{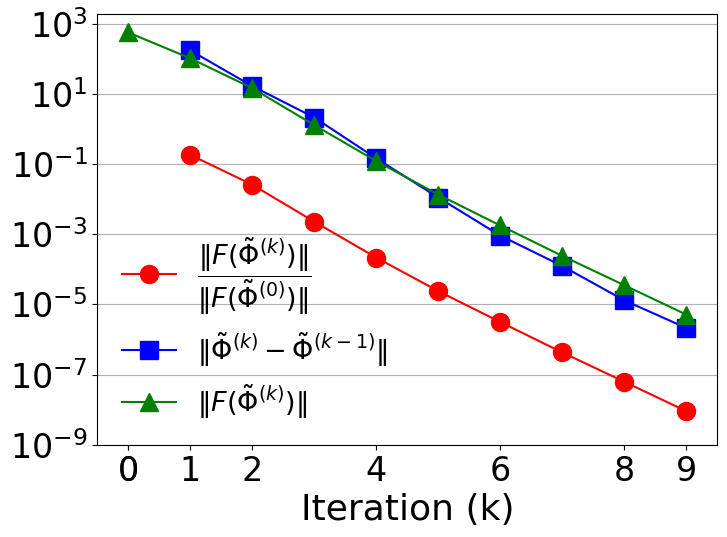}
                \caption{Case of Mesh 4}
               \label{subfig:mesh4-conv}
        \end{subfigure}
        \hspace{0.1cm} 
        \begin{subfigure}[b]{0.315\textwidth}
                \centering
               \includegraphics[width=\textwidth]{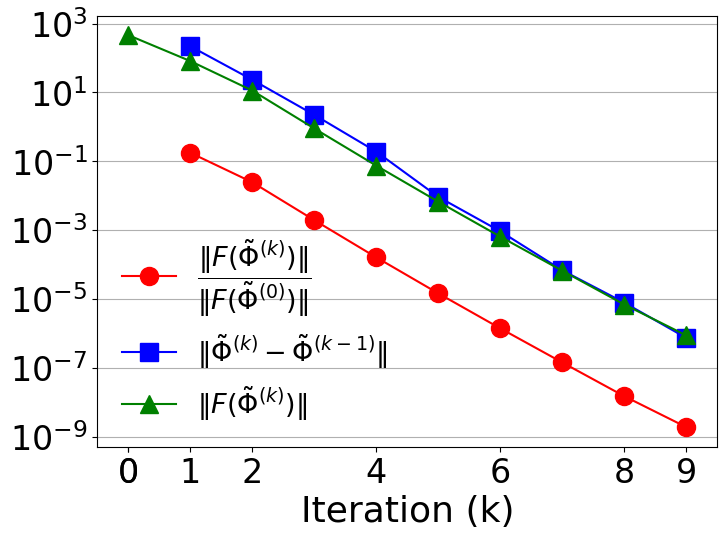}
                \caption{Case of Mesh 5}
               \label{subfig:mesh5-conv}
        \end{subfigure} 
        \hspace{0.1cm}
        \begin{subfigure}[b]{0.315\textwidth}
                \centering
              \includegraphics[width=\textwidth]{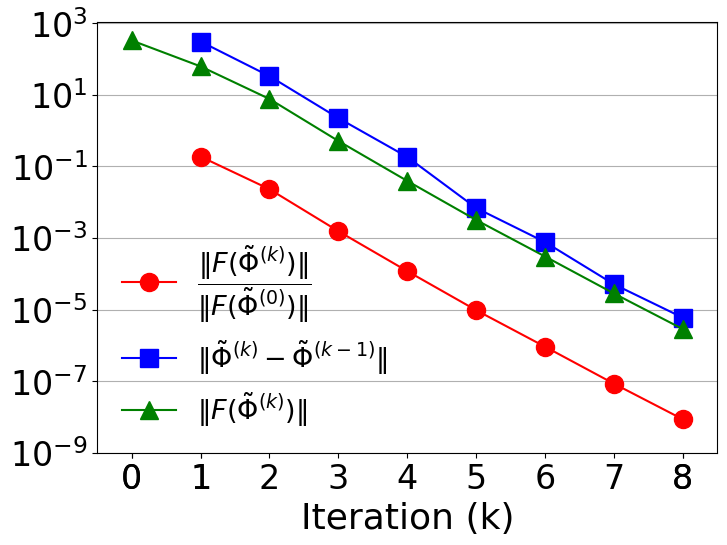}
                \caption{Case of Mesh 6}
               \label{subfig:mesh6-conv}
        \end{subfigure}    
    \caption{Convergence processes of our modified Newton iterative method~\eqref{Newton_scheme2} on Meshes 1 to 6 in the case of protein 1CID.}
    \label{convergence4mesh}
\end{figure}

\begin{table}[t]
\centering
\begin{tabular}{|c|c|c||c|c|}
\hline
Mesh   	&\multicolumn{2}{c||}{Number of iterations} &\multicolumn{2}{c|}{Computer CPU time (seconds)} \\  
        \cline{2-5}
   &  By \eqref{NMPBE-initial1} & By \eqref{NMPBE-initial2} & 
    By  \eqref{NMPBE-initial1} & By \eqref{NMPBE-initial2}  \\
\hline\hline
 \multicolumn{5}{|c|}{Protein (1CID) with 2783 atoms in $\Omega = [-60, 40] \times [-13, 88] \times [-27, 94]$}  \\
    \cline{1-5}
Mesh 1 & 12 & 14 & 1.88 & 1.63 \\ \hline
Mesh 2 & 10 & 11 & 5.31 & 4.67 \\ \hline
Mesh 3 & 11 & 10 & 26.39 & 21.63 \\ \hline
Mesh 4 & 10 & 9 & 51.14 & 43.62 \\ \hline
Mesh 5 & 9 & 9 & 1.45 min. & 1.21 min. \\ \hline
Mesh 6 & 9 & 8 & 4.57 min. & 3.87 min. \\ \hline
\hline
 \multicolumn{5}{|c|}{Protein (4PTI) with 892 atoms in $\Omega = [-23, 51]\times [-13, 56] \times [-30, 44]$}  \\ \cline{1-5} 
Mesh 1 & 11 & 14 & 0.94 & 0.85 \\\hline
Mesh 2 & 11 & 12 & 2.98 & 2.69 \\\hline
Mesh 3 & 10 & 10 & 13.39 & 12.38 \\\hline
Mesh 4 & 10 & 11 & 28.27 & 27.42 \\\hline
Mesh 5 & 9 & 9 & 45.17 & 41.86 \\\hline
Mesh 6 & 8 & 8 & 2.44 min. & 2.07 min. \\\hline
\hline
    \multicolumn{5}{|c|}{Protein (1C4K) with 11439 atoms in $\Omega = [4, 154] \times [-42, 97] \times [-32, 108]$}  \\ \cline{1-5}
\hline
Mesh 1 & 17 & 17 & 17.47 & 5.33 \\ \hline
Mesh 2 & 11 & 11 & 17.06 & 15.46 \\\hline
Mesh 3 & 11 & 10 & 70.68 & 57.67 \\\hline
Mesh 4 & 10 & 10 & 2.55 min. & 2.36 min. \\\hline
Mesh 5 & 9 & 9 & 4.42 min. & 3.63 min. \\\hline
Mesh 6 & 8 & 8 & 14.93 min. & 12.37 min. \\\hline
\end{tabular}
\caption{A comparison of the performance of our modified Newton iterative method \eqref{Newton_scheme2} using the initial iterates $\Phit^{(0)}$ and $\zeta^{(0)}$ generated by \eqref{NMPBE-initial1} and \eqref{NMPBE-initial2}.}
\label{Performance-Newton2}
\end{table}

\begin{figure}[h!]
\centering
\includegraphics[width=0.6\textwidth]{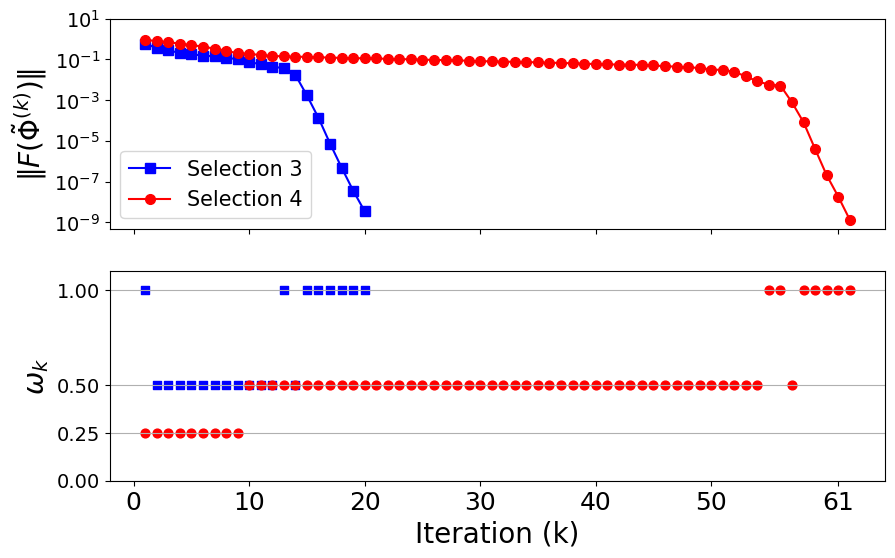}      
\caption{A comparison of the iterative process of the modified Newton method~\eqref{Newton_scheme2} using Selection 3  of initial iteration \eqref{NMPBE-initial3} with that using Selection 4 of \eqref{NMPBE-initial4} for the protein 1C4K case on Mesh 3 in terms of the absolute residual error $\|F(\Phit^{(k)})\|$ and the damping parameter $\omega_k$ of the modified Newton iterative method. 
     }       
\label{fig:line-search}          
\end{figure} 

\subsection{Performance of the NSMPB finite element program package}
Using the meshes described in Table~\ref{table: mesh-data}, we conducted numerical tests on the three proteins 1CID, 4PTI, and 1C4K. The results are presented in Table~\ref{Performance-CPU}, Table~\ref{Performance-Newton2}, and Figure~\ref{convergence4mesh}.

Table~\ref{Performance-CPU} details the CPU time distribution across the five main components of the NSMPB finite element program package, together with the total computation time required from the input of model constants and solver parameters to the output of the results. It also lists the number of mesh vertices, $N_h$, illustrating the growth in mesh size from Mesh 1 to Mesh~6. For each mesh, each nonlinear finite element system has approximately $2N_h$ unknowns because it has two unknown functions: $\Phit$ and $\zeta$ (the convolution of $\Phit$). As shown in Table~\ref{Performance-CPU}, the NSMPB package rapidly generated meshes and efficiently calculated the component functions $G$, $\Psi$, and $\Phit$, the initial iterate $\Phit^{(0)}$, the convolution function $\hat{G}$, and the gradient vectors $\nabla G$ and $\nabla \hat{G}$. 

For example, in the case of protein 1CID using Mesh 5 (with 193,592 vertices), the finite element meshes were generated in approximately 13 seconds, the initial iteration was calculated in 7 seconds, $\Psi$ in 9 seconds, and $\Phit$ in about 1.22 minutes. The total computation time from input to output was approximately 2 minutes. These results demonstrate the high performance of our NSMPB finite element program package in computer CPU time.

Table~\ref{Performance-Newton2} compares the performance of our modified Newton iterative method \eqref{Newton_scheme2} using the two initial iterate selections, Selection 1 of \eqref{NMPBE-initial1} and Selection 2 of \eqref{NMPBE-initial2}, in terms of the number of iterations and CPU time. Notably, both initial iterate selections led to convergence in a similar number of iterations. In these tests, the damping parameter $\omega_k$ was found to be 1 for all $k\geq 0$, enabling our modified Newton iterative method to maintained a fast convergence rate of the Newton iterative method in these tests. From Table~\ref{Performance-Newton2}, it can also been seen that Selection~2 required significantly less CPU time than that of Selection~1. Given its consistent advantage in computational efficiency, Selection~2 has been adopted as the default choice for initial iteration in the NSMPB package.

Figure~\ref{convergence4mesh} displays the convergence processes of our modified Newton iterative method on Meshes 1 to 6, focusing on the protein 1CID case since the results for the other two protein cases are similar. In these tests, initial iterates were generated using Selection~2. In these tests, we evaluated each convergence progress using three types of error: relative residual error $\|F(\Phit^{(k)})\| / \|F(\Phit^{(0)})\|$, difference error $\|\Phit^{(k)} - \Phit^{(k-1)}\|$, and absolute residual error $\|F(\Phit^{(k)})\|$. From Figure~\ref{convergence4mesh}, it can be seen that these three types of errors decreased rapidly, from about 1, $10^2$, and $10^3$ to values close to $10^{-8}$, $10^{-6}$, and $10^{-5}$, respectively, in only 8 to 14 iterations, further demonstrating that our modified Newton iterative method \eqref{Newton_scheme2} has a fast convergence rate.

\subsection{Tests on our damping parameter selection scheme}

We conducted numerical tests in Mesh 3 for the protein 1C4K case using Selections 3 and 4 of the initial iterations given in \eqref{NMPBE-initial3} and \eqref{NMPBE-initial4} to demonstrate the role played by our damping parameter selection scheme in improving the convergence of our modified Newton iterative method. Without using the scheme, the method was found to be divergent in these tests. The test results are reported in Figure~\ref{fig:line-search}.

Figure~\ref{fig:line-search} presents a comparison of the convergence process of our modified Newton iterative method using Selection~3 with that using Selection~4. From Figure~\ref{fig:line-search}, it can be seen that the damping parameter selection scheme, given in Algorithm~1, adjusted the damping parameter $\omega_k$ from 1 to 0.5 or 0.25 to reduce the absolute residual error $\|F(\Phit^{(k)})\|$ per iteration. As soon as the iterates entered a convergence range of Newton's method, we got $\omega_k=1$ for $k\geq n_0$. Here, $n_0$ denotes the number of iterations required to find the convergence region. In these tests, we found $n_0=13$ for Selection~3 and $n_0=56$ for Selection~4 as shown in Figure~\ref{fig:line-search}. When $k\geq n_0$, the modified Newton iterative method quickly converged in 7 iterations for the case of Selection~3 (total 21 iterations) and 6 iterations for the case of Selection~4 (total 62 iterations). These test results demonstrate that our damping parameter selection scheme plays a significant role in improving the convergence and robustness of our modified Newton iterative method.

Recall that we got $n_0=1$ for Selections 1 and 2 in the numerical tests of Section~4.2. Hence, Selections 1 and 2 are much better choices of initial iterations than Selections 3 and 4 in these tests.

\subsection{Tests on finite element solution sequence convergence}

\begin{figure}[h!]
\centering
    \begin{subfigure}{0.32\textwidth}
        \includegraphics[width=\textwidth]{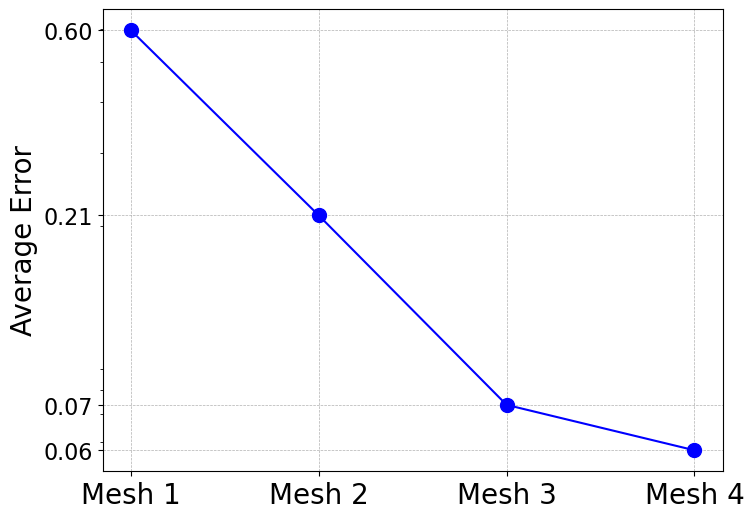}    
    \caption{Case of Protein 4PTI}
\end{subfigure}
\begin{subfigure}{0.32\textwidth}
        \includegraphics[width=\textwidth]{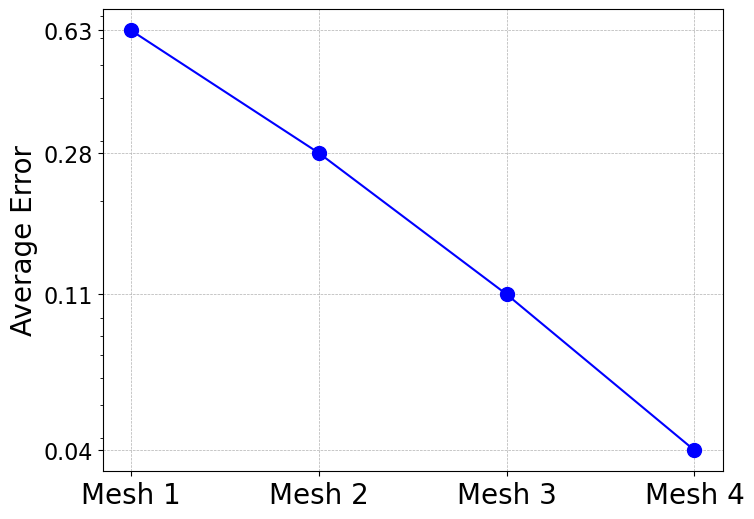}   
    \caption{Case of Protein 1CID}
\end{subfigure}
\begin{subfigure}{0.32\textwidth}
        \includegraphics[width=\textwidth]{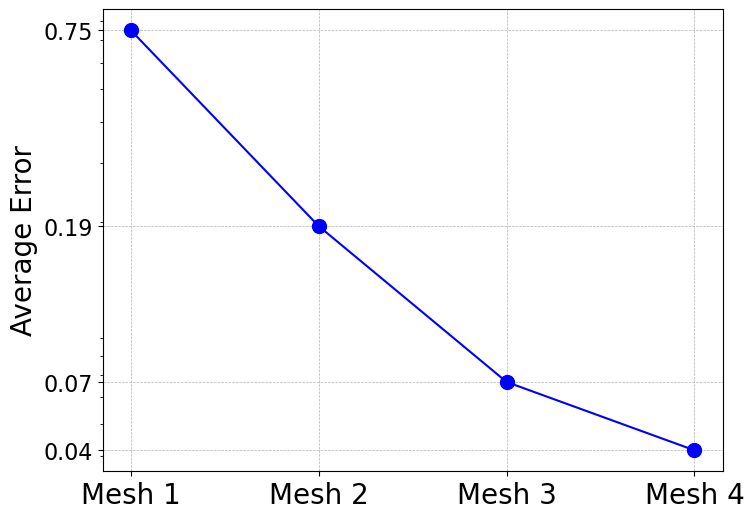}
    \caption{Case of Protein 1C4K}
\end{subfigure}
\caption{Convergence of the NSMPB finite element solutions with increasing mesh size from Mesh 1 to Mesh 4. Here, the average error is defined in \eqref{error_def}, the mesh data for Meshes 1 to 4 are given in Table~\ref{table: mesh-data}, and the solution data produced for Table~\ref{Performance-CPU} are used in the calculation.    
}
\label{figure: convergence}
\end{figure}

\begin{figure}[h!]
\centering
    \begin{subfigure}{0.32\textwidth}
        \includegraphics[width=\textwidth]{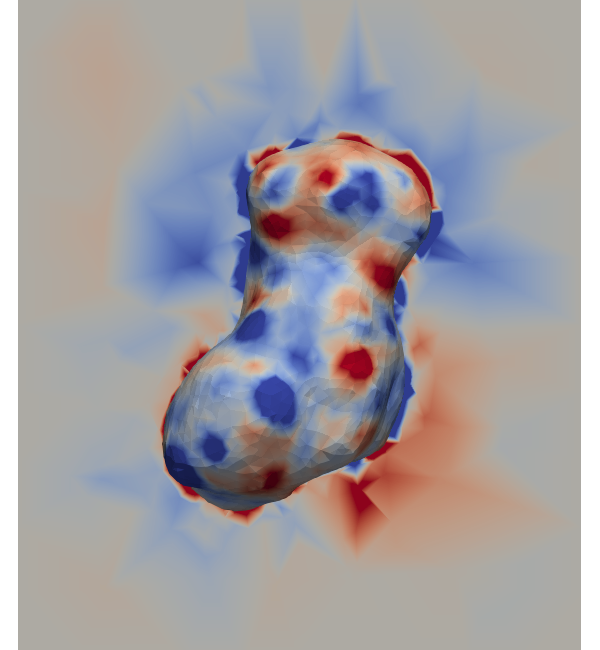}    
    \caption{Case of Mesh 1}
\end{subfigure}
\begin{subfigure}{0.32\textwidth}
        \includegraphics[width=\textwidth]{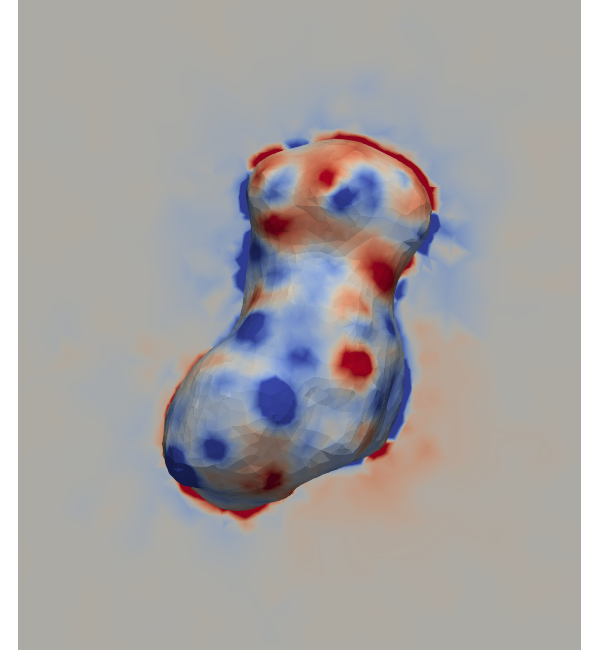}   
    \caption{Case of Mesh 2}
\end{subfigure}
\begin{subfigure}{0.32\textwidth}
        \includegraphics[width=\textwidth]{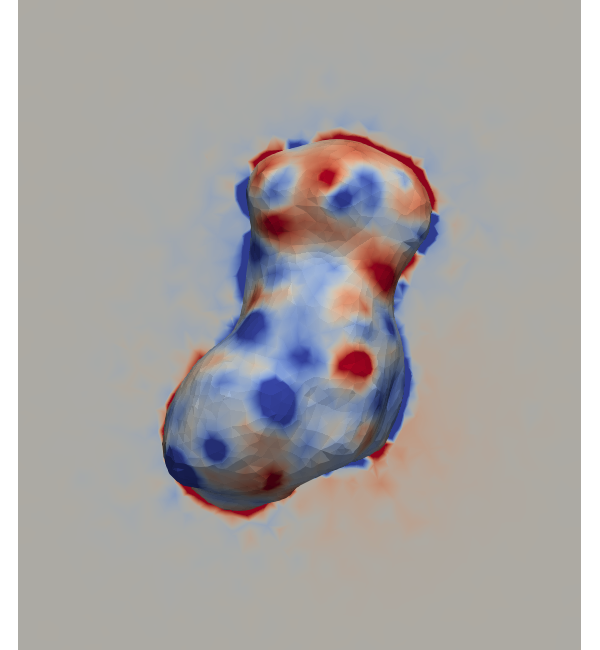}
    \caption{Case of Mesh 3}
\end{subfigure}

\begin{subfigure}{0.32\textwidth}
        \includegraphics[width=\textwidth]{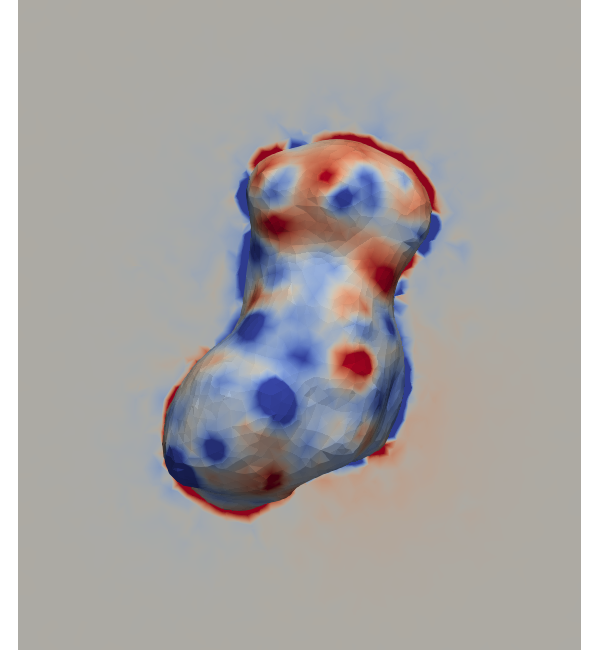}  
    \caption{Case of Mesh 4}
\end{subfigure}
\begin{subfigure}{0.32\textwidth}
        \includegraphics[width=\textwidth]{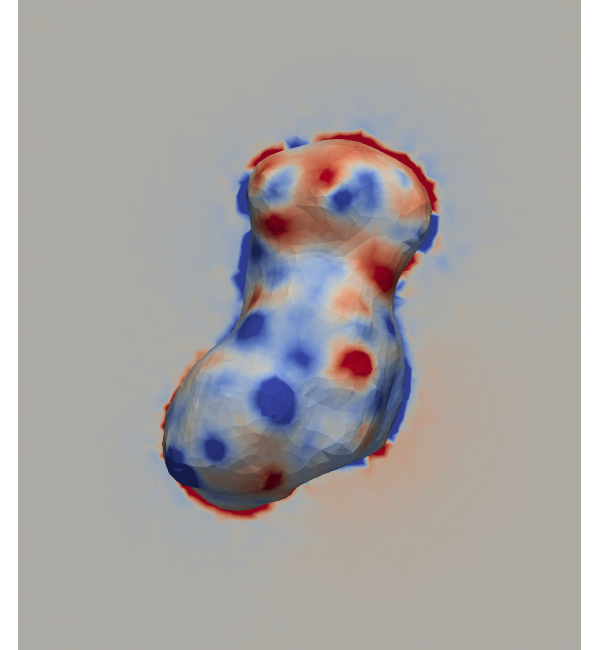}   
    \caption{Case of Mesh 5}
\end{subfigure}
\begin{subfigure}{0.32\textwidth}
        \includegraphics[width=\textwidth]{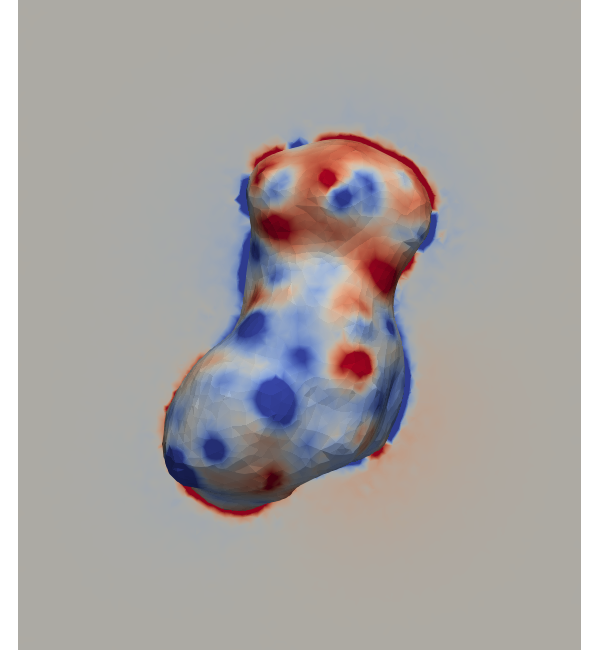}  
    \caption{Case of Mesh 6}
\end{subfigure}
\caption{Color mappings of the six NSMPB finite element solutions, $u_{h_j}$, computed on Meshes 1 through 6 onto the cross-section $x=0$ of the solvent region $D_s$ and a part of the interface between the protein region $D_p$ and $D_s$ in the case of protein 1CID. The color mapping scale for $D_s$ spans from blue corresponding to $u_{h_j} =-1$ to red corresponding to $u_{h_j} =1$, while for $D_p$ from $-5$ to 5.
}
\label{figure: color mapping}
\end{figure}

Let $u_{h_j}$ denote a finite element solution of the NSMPB model generated from the NSMPB package on Mesh $j$. The convergence of a finite element solution sequence, $u_{h_j}$ for $j\geq 1$, can be measured by a sequence of average errors as follows:
\begin{equation}
    \label{error_def}
    \frac{\| u_{h_j} - u^* \|_{L_2(\Omega)}}{N_{h_j}}, \quad j=1, 2, 3, \ldots,
\end{equation}
where $u^*$ denotes a limit function of the sequence, $N_{h_j}$ is the dimension of the linear finite element function space defined on Mesh $j$, and $\| \cdot \|_{L_2(\Omega)}$ is the norm of the function space $L_2(\Omega)$, which is defined by 
\[ \| v \|_{L_2(\Omega)} = \sqrt{\int_{\Omega} |v(\rr)|^2 d\rr} \qquad \forall v\in L_2(\Omega). \]
Because $u^*$ is unknown for the NSMPB model, we substituted the finite element solution $u_{h_6}$ for $u^*$ in the calculation of the average errors for the finite element solutions $u_{h_j}$ for $j=1,2,3$, and 4. 

In fact, in each protein case, $N_{h_6}$ is much larger than $N_{h_j}$ for $j=1$ to 4 as shown in Table~\ref{table: mesh-data}. Moreover, the maximum volume of tetrahedra in Mesh 6 has been limited to 1, resulting in a much higher resolution than Meshes 1 to 4, whose maximum volumes of tetrahedra have been limited to more than 100, 100, 10, and 5, respectively. Thus, $u_{h_6}$ has much higher numerical precision than $u_{h_j}$ for $j=1$ to 4. Hence, it is a good replacement for $u^*$ in the calculation of average errors. We calculated the average errors using the solution data produced for Table~\ref{Performance-CPU}, and reported them in Figure~\ref{figure: convergence}. 

Figure~\ref{figure: convergence} shows that the average errors of the finite element solutions produced with Meshes 1 to 4 consistently decrease, confirming the convergence behavior of a sequence of finite element solutions produced by our NSMPB program package.

To further confirm the convergence of the NSMPB finite element solutions, we plotted the color maps of the six finite element solutions $u_{h_j}$ for $j=1$ to 6 onto the cross section $x=0$ of the solvent region $D_s$ and part of the interface $\Gamma$ between the protein region $D_p$ and $D_s$, which is also part of the boundary of $D_p$. Due to the similarity of these color maps, we only report the case of protein 1CID in Figure~\ref{figure: color mapping}. 

Figure~\ref{figure: color mapping} confirms that the NSMPB finite element solution $u_{h_1}$ in Mesh 1 is significantly less accurate compared to the solution in Mesh 6. The solutions in Meshes 2 through 5 become increasingly similar to $u_{h_6}$ in Mesh 6, indicating improved numerical accuracy with increasing mesh density. Specifically, Mesh 1 is too coarse and Mesh 3 is fine enough to yield a good approximation of the NSMPB solution. Furthermore, the potential values depicted in Figure~\ref{figure: color mapping} align with fundamental electrostatic principles. As expected, regions with negative and positive potential values are predominantly blue and red, indicating areas of strong negative and positive electrostatic potentials. This pattern is particularly noticeable at the interface between the solvent and protein regions across all six meshes, demonstrating the consistency of the NSMPB finite element solutions with basic physical principles.  

\section{Conclusions}

In this paper, we have introduced the nonlocal size modified Poisson-Boltzmann (NSMPB) model to address both ionic size effects and nonlocal dielectric correlations simultaneously, bridging the gap between the capabilities of existing size modified Poisson-Boltzmann (SMPB) and nonlocal modified Poisson-Boltzmann (NMPB) models. To overcome the numerical challenges inherent in the NSMPB model, such as increased nonlinearity, complex nonlocal terms, and solution singularities, we developed a novel solution decomposition and an innovative modified Newton iterative method. The method is further enhanced by an effective damping parameter selection scheme and well-chosen initial iterations. Additionally, we have derived a linear NSMPB model, which offers a computationally inexpensive alternative for protein simulations.

To enable broad applications, we have implemented the NSMPB finite element iterative solver in Python and Fortran as a software package. In this package, we integrate tools for protein structure downloading and processing, interface-fitted tetrahedral mesh generation, parameter assignment, and the output of mesh data files, electrostatic potential functions, and ionic concentration functions in formats compatible with visualization tools such as ParaView. The NSMPB package also seamlessly integrates the Fortran subroutines that we wrote to handle computationally intensive tasks. With these efforts, we have not only simplified the usage of the NSMPB package but also improved the performance of our NSMPB package. 

This work has established a strong foundation for future research and applications. We plan to apply the NSMPB model to the calculation of electrostatic solvation free energies. This will enable us to conduct comparative studies of the NSMPB model with the NMPB model, the SMPB model, and other variants of PB, as well as experimental data produced from biochemical laboratories. We also plan to parallelize the NSMPB package to make it work for large-scale simulations. Furthermore, we will develop nonlocal Poisson-Boltzmann ion channel models and their ion size variants to investigate nonlocal dielectric effects in highly heterogeneous environments. These future works will advance our understanding of biomolecular electrostatics and their critical roles in fundamental biological processes.

\section*{Acknowledgments}
This work was partially supported by the National Science Foundation, USA, through the award number DMS-2153376, and the Simons Foundation, USA, through the research award 711776.


\end{document}